\documentclass[10pt]{amsart}

\usepackage[english]{babel}
\usepackage{pstricks}
\usepackage[applemac]{inputenc}
\usepackage{ae, aeguill}
\usepackage{amsmath,amsthm,amssymb,latexsym,amscd}
\usepackage[mathscr]{eucal}
\usepackage{pstricks}
\usepackage{pst-fill}
\usepackage{xy}
\xyoption{all}
\xyoption{poly}
\usepackage{fancyheadings}
\xyoption{ps}
\UsePSspecials{dvips}

\newif\ifpdf
\ifx\pdfoutput\undefined
\pdffalse 
\else
\pdfoutput=1 
\pdftrue
\fi

\ifpdf
\usepackage[pdftex]{graphicx}
\else
\usepackage{graphicx}
\fi

\newcommand{\ie}{\emph{i.e.}\:}
\newcommand{\cf}{\emph{cf.}\:}

\newcommand{\tors}[1]{\tau^{#1}}

\newcommand{\SU}{\mathrm{SU}(2)}
\newcommand{\SL}{\mathrm{SL}_2(\CC)}
\newcommand{\PSL}{\mathrm{PSL}_2(\CC)}
\newcommand{\SO}{\mathrm{SO}(3)}
\newcommand{\su}{\mathfrak{su}(2)}
\newcommand{\sll}{\mathfrak{sl}_2(\CC)}

\newcommand{\ii}{\mathbf{i}}
\newcommand{\jj}{\mathbf{j}}
\newcommand{\Der}[2]{\mathrm{Der}_{#1} #2}
\newcommand{\Inn}[2]{\mathrm{Inn}_{#1} #2}
\newcommand{\kk}{\mathbf{k}}
\newcommand{\I}{\mathbf{1}}
\newcommand{\ZZ}{\mathbb{Z}}
\newcommand{\CC}{\mathbb{C}}
\newcommand{\IR}{\mathbb{R}}

\newcommand{\Reg}[1]{\mathcal{R}\mathrm{eg}(#1)}
\newcommand{\tangent}[2]{T_{#1} #2}
\newcommand{\TZar}[2]{T^{\mathrm{Zar}}_{#1} #2}
\newcommand{\vol}[1]{\omega^{#1}}
\newcommand{\fox}[2]{\frac{\partial #1}{\partial #2}}

\newcommand{\diff}[3]{\left. \frac{d #1}{d #2}\right\vert_{#3}}
\newcommand{\Diff}[3]{\left. \frac{d}{d #2}{#1}\right\vert_{#3}}

\newcommand{\bord}{\partial}
\newcommand{\lk}{\ell\mathit{k}}
\newcommand{\im}{\mathop{\mathrm{im}}\nolimits}
\newcommand{\coker}{\mathop{\mathrm{coker}}\nolimits}

\begin{document}
\theoremstyle{plain}
\newtheorem{theorem}{Theorem}[section]
\newtheorem*{theorem*}{Theorem}
\newtheorem*{MultLemma}{Multiplicativity Lemma}
\newtheorem*{maintheorem*}{Main Theorem}
\newtheorem{prop}[theorem]{Proposition}
\newtheorem{corollary}[theorem]{Corollary}
\newtheorem{lemma}[theorem]{Lemma}
\newtheorem{claim}[theorem]{Claim}
\newtheorem{fact}[theorem]{Fact}
\theoremstyle{definition}
\newtheorem{definition}{Definition}
\newtheorem{example}{Example}
\newtheorem*{example*}{Example}
\newtheorem{notation}{Notation}
\newtheorem*{notation*}{Notation}
\newtheorem*{convention*}{Convention}
\theoremstyle{remark}
\newtheorem{remark}{Remark}
\newtheorem{question}{Question}
\pagestyle{fancy}
\lhead{}\chead{}
\chead{To appear in \emph{Annales de l'Institut Joseph Fourier}}
%
\date{\today}
\title{Non abelian Reidemeister torsion and volume form on the $\SU$-representation space of knot groups}
\author{J\'er\^ome Dubois -- \emph{\textrm{\today}}}
\address{Section de Math\'ematiques \\ Universit\'e de Gen\`eve CP 64, 2--4 Rue du Lièvre \\ CH 1211 Genève 4 (Suisse)}
\email{Jerome.Dubois@math.unige.ch}
\ifpdf
\DeclareGraphicsExtensions{.pdf, .jpg, .tif}
\else
\DeclareGraphicsExtensions{.eps, .jpg}
\fi
\begin{abstract}For a knot $K$ in $S^3$ and a {regular} representation $\rho$ of its group $G_K$ in $\SU$ we construct a non abelian Reidemeister torsion form on the first twisted cohomology group of the knot exterior. 
This non abelian Reidemeister torsion form provides a volume form on the $\SU$-representation space of $G_K$ (see Section~\ref{Torsion}). In another way, we construct according to Casson---or more precisely taking into account Lin's~\cite{Lin:1992} and Heusener's~\cite{Heu:2003} further works---a volume form on the $\SU$-representation space of $G_K$ (see Section~\ref{Vol}). Next, we compare these two apparently different points of view---the first by means of the Reidemeister torsion and the second defined ``à la Casson"---and finally prove that they produce the {same} topological knot invariant (see Section~\ref{Comparison}). 
\end{abstract}
\maketitle

\markboth{Jérôme Dubois}{Non abelian Reidemeister torsion and volume form}
\pagestyle{myheadings}

\section*{Introduction: Motivation and Main ideas}

	The space of conjugacy classes of irreducible $\SU$-representations of a knot group is a real semi-algebraic set and its {regular} part is a $1$-dimensional ma\-ni\-fold (for more details see \emph{infra}).  The aim of this paper is to compare two \emph{a priori} different constructions on this $1$-dimensional manifold: one by means of Rei\-de\-meis\-ter torsion and another using Casson's original construction. Firstly, we associate to any regular $\SU$-representation $\rho$ of the knot group a non abelian Reidemeister torsion form on the first cohomology group of the knot exterior with coefficients in the adjoint representation associated to $\rho$. Secondly, we construct ``à la Casson" a volume form on the $1$-dimensional manifold consisting of conjugacy classes of regular representations, which appears to be a knot invariant. Finally, we prove that these \emph{a priori} two different points of view are indeed {equivalent} (see~Theorem~\ref{T:Volume=Torsion}). \\

	The {Reidemeister torsion} was introduced for the first time in 1935 by K. Rei\-de\-meis\-ter in his work~\cite{Reidemeister:1935} on the combinatorial classification of $3$-dimensional lens spaces. The torsion is a combinatorial invariant which is not a homotopy invariant. It has became a fundamental tool in low-dimensional topology, see for example Turaev's monograph~\cite{Turaev:2002}. Informally speaking, the Reidemeister torsion is a graded version of the determinant in the same way that the Euler characteristic is a graded version of the dimension.

	In 1985, A. Casson constructed an integer valued invariant of integral homology $3$-spheres which appeared extremely useful. 
	The original definition of the Casson invariant is based on $\SU$-representation spaces. Informally speaking, the Casson invariant of the homology sphere $M$ counts algebraically the number of conjugacy classes of irreducible $\SU$-representations of $\pi_1(M)$ in the same sense that the Lefschetz number of a map counts the number of fixed points. 
	
	In 1992, X.-S. Lin used Casson's construction to define an integer knot invariant. He indirectly proved that this invariant is equal to half the signature of knots (see~\cite{Lin:1992}). At first sight, the equality between these two apparently different quantities seems mysterious. In 2003, M. Heusener explained Lin's result using an orientation on the $\SU$-representation space of knot groups. More precisely, given a knot $K \subset S^3$ let $M_K$ denote its exterior and let $G_K = \pi_1(M_K)$ denote its group. In general the representation space of $G_K$ in $\SU$ has singularities; to avoid this difficulty Heusener and Klassen introduced in~\cite{HeuKlassen:1997} the notion of {regular} representation. 
	An irreducible $\SU$-representation $\rho$ of $G_K$ is called {regular} if $\dim H^1_\rho(M_K) = 1$, where $H^*_\rho(M_K)$ denotes the $({Ad\circ \rho})$-twis\-ted cohomology of $M_K$. Let $\Reg{K}$ denote the set of conjugacy classes of regular representations of $G_K$. Heusener proved that $\Reg{K}$ is a canonically {oriented} $1$-dimensional manifold  (see~\cite[Section~1]{Heu:2003}).\\

	In this article we investigate a {volume form} on $\Reg{K}$ which is an invariant of $K$. In particular, we explain an intimate connection between this volume form and a non abelian Reidemeister torsion form on $H^1_\rho(M_K)$. 
	Observe that our torsion is different from the one previously studied by E. Witten for surfaces in~\cite{Witten:1991}. This difference is precisely due to the non-triviality of the second twisted cohomology group $H^2_\rho(M_K)$ of $M_K$. The definition of the Reidemeister torsion form for the exterior of $K$  needs a reference basis for $H^2_\rho(M_K)$. Our construction of this basis is inspired by the work~\cite[Corollary 3.23]{Porti:1997} of J. Porti.  
	Let us describe the main ideas of the construction. For a {regular} representation $\rho$ of $G_K$ in $\SU$, the twisted cohomology groups of $M_K$ satisfy: $H^0_\rho(M_K) = 0$, $H^1_\rho(M_K)\cong \IR$ and $H^2_\rho(M_K)\cong \IR$. We produce a {distinguished basis} for $H^2_\rho(M_K)$ using the fundamental class of $\bord M_K$ and the isomorphisms $H^2_\rho(M_K) \cong H^2_\rho(\bord M_K) \cong H^2(\bord M_K; \ZZ) \otimes \IR$. With this reference basis, we construct a {linear form} $\tau^K_\rho$ on the $1$-dimensional vector space $H^1_\rho(M_K)$. This form is called the {Reidemeister torsion form}  associated to $K$ and $\rho$. It is a combinatorial invariant and maps any generator of $H^1_\rho(M_K)$ to the sign-determined $({Ad\circ \rho})$-twisted Reidemeister torsion of $M_K$ corresponding to this choice of bases for $H^*_\rho(M_K)$ and to the canonical cohomology orientation of $M_K$. 

	In another way, we construct according to Casson a canonical volume form on $\Reg{K}$. This volume form is denoted $\omega^K$. The construction of $\omega^K$ uses a {plat} presentation of the knot $K$ to obtain a splitting of its exterior. This particular decomposition replaces the Heegaard splitting used in the original construction of the Casson invariant for homology $3$-spheres (see~\cite{GM:1992}). 
	
	Let $\widehat{R}(M_K)$ denote the space of conjugacy classes of irreducible $\SU$-repre\-sen\-ta\-tions of $G_K$ and assume that $\rho$ is a regular representation. We show that the tangent space to $\widehat{R}(M_K)$ at the conjugacy class of $\rho$ is isomorphic to $H^1_\rho(M_K)$. This isomorphism allows us to formulate the main result of the present work.	
\begin{maintheorem*}[Theorem~\ref{T:Volume=Torsion}]
If $\rho: G_K \to \SU$ is regular,  then ${\tau}^K_{\rho} = \vol{K}_{[\rho]}.$
\end{maintheorem*}

	The present paper is organized as follows. Sections~\ref{RepSp} to~\ref{TwistDer} are reviews on $\SU$-representation spaces, sign-determined Reidemeister torsion and twisted derivations respectively. In Section~\ref{TwistH}, we discuss the canonical cohomology orientation and the twisted cohomology of knot exteriors and also the notion of regular representations. In Section~\ref{Torsion}, we investigate the non abelian Reidemeister torsion form. Section~\ref{Vol} deals with the Casson-type construction of the natural volume form on the $\SU$-re\-pre\-sen\-ta\-tion space of knot groups. In Section~\ref{Comparison}, we explain how to interpret this volume form in terms of the Reidemeister torsion form.  The rest of the article is devoted to the proof of the Main Theorem (Theorem~\ref{T:Volume=Torsion}), see Sections~\ref{S:Volume=Torsion1} and~\ref{S:Volume=Torsion2}.
	
	This paper is part of the author  Ph.D. thesis~\cite{JDTHESE}.

\section{Representation spaces}\label{RepSp}

	Here we collect some well-known results about $\SU$-representation spaces and introduce the notation used throughout this paper.

\subsection{Some notation}
	In all this article ``\emph{connected}" means ``\emph{arc-connected}".
	
	The fundamental group $\pi_1(X)$ of a connected topological space $X$ is consider without specifying a base point since all the constructions we do are invariant by conjugation, see~\cite[p. 9]{Porti:1997} for details.
	
	The Lie group $\SU$ acts on its Lie algebra $\su$ via the adjoint representation $Ad_A : \su \to \su$ defined by $Ad_A(x) = AxA^{-1}$, where $A \in \SU$. As a manifold $\SU$ is identified with the $3$-sphere $S^3$ and we identify the $2$-sphere $S^2$ with the set of zero-trace matrices of $\SU$: $S^2 = \{A \in \SU \;|\; \mathrm{Tr}(A) = 0\}$. The Lie algebra $\su$ is identified with the {pure quaternions} i.e. with the quaternions of the form $q = a \ii + b \jj + c \kk$. Recall that for all $A \in \SU$ there exist $\theta \in [0, \pi]$ and $P \in S^2$ such that $A = \cos(\theta) + \sin(\theta)P$. Moreover, the pair $(\theta, P)$ is {unique} if and only if $A \ne \pm \I$. Note that $Ad_A$ is the rotation of angle $2 \theta$ which fixes $P$. The Lie algebra $\su$ is equipped with the usual scalar product defined by $\langle x, y \rangle = -1/2 \cdot \mathrm{Tr}(xy)$. This scalar product coincides with the Killing form of $\su$ multiplied by $-2$. We think of $\SO$ as the base space of the 2-fold covering $Ad : \SU \to \SO$ given by $A \mapsto Ad_A$.

\subsection{Representation spaces}
	Given a finitely generated group $G$ we let $R(G) = \mathrm{Hom}(G; \SU)$ denote the space of $\SU$-representations of $G$. This space is endowed with the {compact-open topology}. Here $G$ is assumed to have the discrete topology and $\SU$ the usual one. A representation $\rho$ of $G$ in $\SU$ is called {abelian} or {reducible} (resp. {central}) if its image $\rho(G)$ is an abelian subgroup of $\SU$ (resp. is contained in the center $\{\pm \I\}$ of $\SU$). A representation is called {irreducible} if it is not abelian. We let $\widetilde{R}(G)$ (resp. $C(G)$) denote the subspace of irreducible (resp. central) representations.
	
	The compact Lie group $\SU$ acts on $R(G)$ by conjugation. We write $[\rho]$ for the conjugacy class of the representation $\rho \in R(G)$ and we let $\SU(\rho)$ denote its orbit. The action by conjugation factors through $\SO=\SU/\{\pm \I\}$ as a {free} action on the open  subspace $\widetilde{R}(G)$ and we set $\widehat{R}(G)=\widetilde{R}(G)/\SO$. In this way, we can think of $\widehat{R}(G)$ as the base space of a principal $\SO$-bundle with total space $\widetilde{R}(G)$, see~\cite[Section~3.A]{GM:1992}.  Furthermore, we know that $\widehat{R}(G)$ admits the structure of a {real affine semi-algebraic set} (see~\cite[Section 2]{Heu:2003}). Here a subset of $\IR^m$ is called a semi-algebraic set if it is a finite union of finite intersections of sets defined by polynomial equalities and inequalities.
	
\begin{notation*}
	For a connected CW-complex $W$ we write $R(W) = R(\pi_1(W))$, $\widetilde{R}(W) = \widetilde{R}(\pi_1(W))$, $\widehat{R}(W) = \widehat{R}(\pi_1(W))$ etc. 
\end{notation*}

	For a knot $K$ in $S^3$ let $M_{K} = S^{3}\setminus N(K)$ denote its {exterior} and let $G_{K} = \pi_{1}(M_{K})$ denote its {group}. Here $N(K)$ is an open tubular neighbourhood of $K$. Recall that $M_{K}$ is a compact $3$-manifold bounded by the $2$-torus $\bord M_K$. 

	Let $G_K'$ be the subgroup generated by the commutators of $G_K$. The abelianization $G_{K}/G_{K}' \cong H_1(M_K; \ZZ)$ is generated by the meridian $m$ of $K$. So each {abelian} representation of $G_{K}$ in $\SU$ is conjugate to one and only one of the $\varphi_\theta : G_{K} \to \SU$ defined by $\varphi_\theta(m) = \cos(\theta)  + \sin(\theta) \ii$, with $0 \leqslant \theta \leqslant \pi$. 

  \begin{notation*}
	Let $\rho \in R(G)$ and let $g \in G$ be such that $\rho(g) \ne \pm \mathbf{1}$. There exists a unique pair $\left({\theta^\rho(g), P^\rho(g)}\right) \in (0, \pi) \times S^2$ such that $\rho(g) = \cos(\theta^\rho(g)) + \sin(\theta^\rho(g))P^\rho(g).$
\end{notation*}

\section{Review on Reidemeister torsion}
\label{SS:RTorsion}

	This section reviews the definitions of a volume form, of the notion of compa\-ti\-bi\-li\-ty and of the sign-determined Reidemeister torsion of a CW-complex to set up the conventions which will be used. For more details, we refer to Milnor's survey ~\cite{Milnor:1966} and to Turaev's monographs~\cite{Turaev:2000}~\&~\cite{Turaev:2002}.

\subsection{Volume forms and compatibility}
\label{S:Vol}
	 
	Let $E$ be a $n$-dimensional real vector space. Let $E^* = \mathrm{Hom}_{\IR}(E, \IR)$ denote its dual space. A {volume form} $v$ on $E$ is a generator of the $n$th exterior power $\bigwedge^n E^*$. 
	Let $E', E''$ be two real finite dimensional vector spaces and let $v', v''$ be volume forms on $E', E''$ respectively. The direct sum $E' \oplus E''$ has a canonical volume form denoted $v' \wedge v''$.
	
	Consider now a short exact sequence $\xymatrix@1@-.7pc{0 \ar[r] & E' \ar[r]^-i & E \ar[r]^-j & E'' \ar[r] & 0}$ of finite dimensional real vector spaces. Let $v'$, $v$ and $v''$ be volume forms on $E'$, $E$ and $E''$ respectively; let $s$ denote a section of $j$, so that $i \oplus s : E' \oplus E'' \to E$ is an isomorphism. We say that the previous three volume forms are \emph{compatible} with each other if $v' \wedge v'' = (i \oplus s)^*(v).$
	It is easy to verify that the notion of compati\-bi\-li\-ty is independent of the chosen section $s$. Furthermore, if any two of the vector spaces $E'$, $E$ and $E''$ are endowed with a volume form, then the third is endowed with a unique well-defined volume form which is compatible with the two others. In particular, if $v, v''$ are volume forms on $E, E''$ respectively, we write $v' = v/v''$ the unique compatible volume form on $E'$ to indicate its dependence on $v$ and $v''$ (\cf~\cite[Section 1]{Milnor:1962}). Compatibility will be used to set up the Multiplicativity Lemma (see Subsection~\ref{SS:Multiplicativite}) as well as in order to build up  ``new" volume forms (see Subsection~\ref{SS:Construction}).

\subsection{Algebraic preliminaries}\label{AlgPre}
	
	Keep the notation of the previous subsection.
	
	For two ordered bases $\mathbf{a} = \{a_1, \ldots, a_n\}$ and $\mathbf{b} = \{b_1, \ldots, b_n\}$
 of $E$ we write $[\mathbf{a}/\mathbf{b}] = \det (p_{i j})_{i, j}$, where $a_i = \sum_{j=1}^n p_{i j} b_j$, for all $i$. The bases $\mathbf{a}$ and $\mathbf{b}$ are {equivalent} if $[\mathbf{a}/\mathbf{b}]= +1$. Two equivalent bases define the same volume form on $E$.
 	
	Let $C_* = ( \xymatrix@1@-.5pc{0 \ar[r] & C_n \ar[r]^-{d_n} & C_{n-1} \ar[r]^-{d_{n-1}} & \cdots \ar[r]^-{d_1} & C_0 \ar[r] & 0})$ be a chain complex of finite dimensional vector spaces over $\IR$. For each $i$, consider $B_i = \im(d_{i+1} : C_{i+1} \to C_i)$, $Z_i = \ker(d_{i} : C_i \to C_{i-1})$ and the homology group $H_i = Z_i / B_i$.	
	 We suppose that $C_i$ (resp. $H_i$) is endowed with a {reference basis} $\mathbf{c}^i$ (resp. $\mathbf{h}^i$) for each $i$.  In this case, we say that $C_*$ is {based} and {homology based}. The Reidemeister torsion of the based and homology based chain complex $C_*$, with reference bases $\mathbf{c}^*$ and $\mathbf{h}^*$, is defined as follows. Let $\mathbf{b}^i$ denote a basis of $B_i$. We have the two short exact sequences:
\begin{equation}\label{SE1}
\xymatrix@1@-.5pc{0 \ar[r] & Z_i \ar[r] & C_i \ar[r] & B_{i-1} \ar[r] & 0,}
\end{equation}
\begin{equation}\label{SE2}
\xymatrix@1@-.5pc{0 \ar[r] & B_i \ar[r] & Z_i \ar[r] & H_i \ar[r] & 0.}
\end{equation}
Choose a lift $\widetilde{\mathbf{b}}^i$ of $\mathbf{b}^i$ in $C_i$ and a lift $\widetilde{\mathbf{h}}^i$ of $\mathbf{h}^i$ in $Z_i$. Many choices of such lifts are of course possible. Using the exactness of~(\ref{SE1}) and~(\ref{SE2}) the sequences of vectors $\mathbf{b}^i$, $\widetilde{\mathbf{h}}^i$ and $\widetilde{\mathbf{b}}^{i-1}$ combine to yield a new basis $\mathbf{b}^i  \widetilde{\mathbf{h}}^i  \widetilde{\mathbf{b}}^{i-1}$ of $C_i$. 
	With this notation, the {Reidemeister torsion} of $C_*$, with reference bases $\mathbf{c}^*$ and $\mathbf{h}^*$, is (see~\cite[Definition 3.1]{Turaev:2000}): 
\begin{equation}
\label{Def:RTorsion}
\mathrm{tor}(C_*, \mathbf{c}^*, \mathbf{h}^*) = \prod_{i=0}^n [\mathbf{b}^i  \widetilde{\mathbf{h}}^i \widetilde{\mathbf{b}}^{i-1}/\mathbf{c}^i]^{(-1)^{i+1}} \in \IR^*.
\end{equation}
 The torsion $\mathrm{tor}(C_*, \mathbf{c}^*, \mathbf{h}^*)$ does not depend on the choice of the bases $\mathbf{b}^i$ of $B_i$ nor on the lifts $\widetilde{\mathbf{b}}^i$ and $\widetilde{\mathbf{h}}^i$.

	We say that the based chain complex $C_*$ is {acyclic} if $H_i(C_*)$ vanishes for all $i$. In this case, if $\Gamma_* = (C_{*} \to 0)$ denotes the left-shift chain complex obtained from $C_{*}$, \ie $\Gamma_0 = 0$ and $\Gamma_i = C_{i-1}$ for $i\geqslant 1$, then 
\begin{equation}\label{EQ:Decalage}
\mathrm{tor}(\Gamma_*) = (\mathrm{tor}(C_*))^{-1}.
\end{equation}

	The Reidemeister torsion of $C_*$ does only depend on the equivalence classes of the reference bases $\mathbf{c}^i$ and $\mathbf{h}^i$. More precisely, if $\mathbf{c'}^i$ is a different basis of $C_i$ and $\mathbf{h'}^i$ a different one of $H_i$, then we have the \emph{basis change formula}
\begin{equation}
\label{EQ:changementdebase}
\frac{\mathrm{tor}(C_*, \mathbf{c'}^*, \mathbf{h'}^*)}{\mathrm{tor}(C_*, \mathbf{c}^*,  \mathbf{h}^*)} = \prod_{i=0}^n \left( \frac{[{\mathbf{c}'}^i/\mathbf{c}^i]}{[{\mathbf{h}'}^i/\mathbf{h}^i]}\right)^{(-1)^i}.
\end{equation}

\begin{convention*}
	Suppose that the chain complex $C_*$ is not based and not homology based but is such that each $C_i$ and each $H_i$ is endowed with distinguished volume forms. The Reidemeister torsion of such $C_*$ is computed---according to formula~(\ref{EQ:changementdebase})---with respect to any reference bases of $C_i$ and $H_i$ which have volume one.
\end{convention*}

\subsection{The Reidemeister torsion of a CW-complex}

	If formula~(\ref{Def:RTorsion}) is used to define the Reidemeister torsion of a CW-complex, then we will fall into the well-known ``{up-to-sign ambiguity}" of the Reidemeister torsion. To solve this problem V.~Turaev has introduced a {sign-determined} Reidemeister torsion. 

\subsubsection*{\textbf{The sign-determined torsion}}

Set 
\[
\alpha_i(C_*) = \sum_{k=0}^i \dim C_k \in \ZZ/2\ZZ,  \;
    \beta_i(C_*)  = \sum_{k=0}^i \dim H_k \in \ZZ/2\ZZ, 
\]
\[
|C_*| = \sum_{k\geqslant 0} \alpha_k(C_*) \beta_k(C_*) \in \ZZ/2\ZZ.
\]
The {sign-determined Reidemeister torsion} of $C_*$ is the ``sign-corrected" torsion
\begin{equation}
\label{EQ:TorsionRaff}
	\mathrm{Tor}(C_*, \mathbf{c}^*, \mathbf{h}^*) = (-1)^{|C_*|}  \, \mathrm{tor}(C_*, \mathbf{c}^*, \mathbf{h}^*) \in \IR^*,
\end{equation}
see~\cite[Section 3.1]{Turaev:1986} or~\cite[formula (1.a)]{Turaev:2002}.

\subsubsection*{\textbf{The Reidemeister torsion of a CW-complex}}

	Let $W$ be a finite CW-complex; consider a representation $\rho : \pi_1(W) \to \SU$. The universal covering $\widetilde{W}$ of $W$ is endowed with the induced CW-complex structure and the group $\pi_1(W)$ acts on $\widetilde{W}$ by the covering transformations. This action turns $C_*(\widetilde{W}; \ZZ)$ into a chain complex of left 
$\ZZ[\pi_1(W)]$-modules. Further observe that the Lie algebra $\su$ can be viewed as a left $\ZZ[\pi_1(W)]$-module via the representation $Ad \circ \rho$. The {$(Ad \circ \rho)$-twisted cochain complex} of $W$ is
\[
C^*(W; Ad \circ \rho) = \mathrm{Hom}_{\pi_1(X)}(C_*(\widetilde{W}; \ZZ), \su).
\]
The cochain complex $C^*(W; Ad \circ \rho)$ computes the {$(Ad \circ \rho)$-twisted cohomology} of $W$. This cohomology is denoted $H^*_\rho(W)$. 
When $H^*_\rho(W) = 0$ we say that $\rho$ is {acyclic}.

	Let $\{e^{(i)}_1, \ldots, e^{(i)}_{n_i}\}$ denote the set of $i$-dimensional cells of $W$. Choose a lift $\tilde{e}^{(i)}_j$  of the cell $e^{(i)}_j$ in $\widetilde{W}$ and choose an arbitrary order and an arbitrary orientation for the cells $\tilde{e}^{(i)}_j$. Thus, for each $i$, $\mathbf{c}^{i} = \{ \tilde{e}^{(i)}_1, \ldots, \tilde{e}^{(i)}_{n_i} \}$ is a $\ZZ[\pi_1(W)]$-basis of $C_i(\widetilde{W}; \ZZ)$ and we consider the corresponding ``dual" basis over $\mathbb{R}$ $$\mathbf{c}^{i}_{\su} = \left\{ {\tilde{e}^{(i)}_{1, \ii}, \tilde{e}^{(i)}_{1, \jj}, \tilde{e}^{(i)}_{1, \kk}, \ldots, \tilde{e}^{(i)}_{n_i, \ii}, \tilde{e}^{(i)}_{n_i, \jj}, \tilde{e}^{(i)}_{n_i, \kk}}\right\}$$ of $C^i(W; Ad \circ \rho) = \mathrm{Hom}_{\pi_1(X)}(C_i(\widetilde{W}; \ZZ), \su)$. 
	
If $\mathbf{h}^{i}$ is a basis of $H^i_\rho(W)$ then $\mathrm{Tor}(C^*(W; Ad\circ \rho), \mathbf{c}^*_{\su}, \mathbf{h}^{*}) \in \IR^*$ is well-defined.

	The cells $\{\tilde{e}^{(i)}_j\}_{0 \leqslant i \leqslant \dim W, 1 \leqslant j \leqslant n_i}$ are in one-to-one correspondence with the cells of $W$ and their order and orientation induce an order and an orientation for the cells $\{e^{(i)}_j\}_{0 \leqslant i \leqslant \dim W, 1 \leqslant j \leqslant n_i}$. We thus produce a basis over $\IR$ for $C^*(W; \IR)$ which is denoted $c^*$. 
	
	Choose a \emph{cohomology orientation} of $W$ \ie an orientation of the real vector space $H^*(W; \IR) = \bigoplus_{i\geqslant 0} H^i(W; \IR)$; let $\mathfrak{o}$ denote such an orientation. Provide each vector space $H^i(W; \IR)$ with a reference basis $h^i$ such that the basis $h^* = \{h^0, \ldots, h^{\dim W}\}$ of $H^*(W; \IR)$ is {positively oriented} with respect to the cohomology orientation $\mathfrak{o}$. Compute the sign-determined Reidemeister torsion $\mathrm{Tor}(C^*(W; \IR), c^*, h^{*}) \in \IR^*$  of the resulting based and cohomology based chain complex $C^*(W; \IR)$ and consider its sign $\tau_0 = \mathrm{sgn}\left(\mathrm{Tor}(C^*(W; \IR), c^*, h^{*})\right) \in \{\pm 1\}$. Further observe that $\tau_0$ does not depend on the choice of the positively oriented basis $h^*$. 
	
	The {sign-determined Reidemeister torsion} of the cohomology oriented CW-com\-plex $W$ twisted by the representation $Ad \circ \rho$ is the product
\[
\mathrm{TOR}(W; Ad\circ \rho, \mathbf{h}^{*}, \mathfrak{o}) = \tau_0 \cdot \mathrm{Tor}(C^*(W; Ad\circ \rho), \mathbf{c}^*_{\su}, \mathbf{h}^{*}) \in \IR^*.
\]
	The torsion $\mathrm{TOR}(W; Ad\circ \rho, \mathbf{h}^{*}, \mathfrak{o})$  is called the \emph{$(Ad \circ \rho)$-twisted Reidemeister torsion} of $W$. It is well-defined. It does not depend on the choice of the lifts $\tilde{e}^{(i)}_j$ nor on the order and orientation of the cells (because they appear twice). Finally, it just depends on the conjugacy class of $\rho$.
	
	One can prove that $\mathrm{TOR}$ is invariant under cellular subdivision, homeomorphism class and simple homotopy type. In fact, it is precisely the sign $(-1)^{|C_*|}$ in~(\ref{EQ:TorsionRaff}) which ensures all these particularly important properties of invariance (see~\cite[Chapter 2]{JDTHESE} for detailed proofs).

\section{Twisted derivations and Tangent bundle to representation spaces}\label{TwistDer}

	This section reviews the useful concept of {twisted derivations}. 

\subsection{Twisted derivations}
	
	Let $G$ be a finitely generated group and consider a representation $\rho : G \to \SU$. 
	
	An ${(Ad\circ \rho)}$-{twisted derivation} is a map $d : G \to \su$ satisfying the cocycle condition: $d(g_1g_2) = d(g_1) + Ad_{\rho(g_1)} d(g_2),$ for all $g_1, g_2 \in G$. We let $\Der{\rho}{(G)}$ denote the set of ${(Ad\circ \rho)}$-twisted derivations. Among twisted derivations we distinguish the inner ones. A map $\delta : G \to \su$ is an {inner derivation} if there exists $a \in \su$ such that $\delta(g) = a - Ad_{\rho(g)} a,$ for all $g\in G$. We let $\Inn{\rho}{(G)}$ denote the set of inner derivations. Observe that $\Inn{\rho}{(G)} \cong \su$ for any {irreducible} representation $\rho$.

	Recall that $$Z^1_\rho(G) \cong \Der{\rho}{(G)}, \; B^1_\rho(G) \cong \Inn{\rho}{(G)}, \; H^1_\rho(G) \cong \Der{\rho}{(G)}/\Inn{\rho}{(G)},$$ and 
$$H^0_\rho(G) = \su^{Ad \circ \rho(G)} = \{v \in \su \; |\; v = Ad_{\rho(g)}v, \; \forall g \in G\}.$$ 
For each irreducible representation $\rho$ of $G$ we thus have the short exact sequence
\begin{equation}\label{ExSH}
\xymatrix@1@-.6pc{0 \ar[r] & \su \ar[r] & \Der{\rho}{(G)} \ar[r] & H^1_\rho(G) \ar[r] & 0.}
\end{equation}

\subsection{Derivations and tangent spaces}
\label{DerTan}
	We begin by a digression on $\SL$-re\-pre\-sen\-tations. Let $R(G, \SL) = \mathrm{Hom}(G; \SL)$ denote the space of $\SL$-re\-pre\-sentations of $G$ endowed with the compact-open topology. To each re\-pre\-sen\-ta\-tion $\rho \in R(G, \SL)$ we associate its {character} $\chi_\rho : G \to \CC$ defined by $\chi_\rho(g) = \mathrm{Tr}(\rho(g))$ for all $g \in G$. The set of {characters} of $G$ is called the {re\-pre\-sen\-ta\-tion variety} and is denoted $X(G)$.  In a way $X(G)$ is the ``algebraic quotient" of $R(G, \SL)$ by the action by conjugation of $\PSL$ because the naive quotient $R(G, \SL)/\PSL$ is not Hausdorff in general. It is well known that $R(G, \SL)$ and $X(G)$ have the structure of {complex algebraic affine sets} (see~\cite{CS:1983}).

	In~\cite{Weil:1964}, A.~Weil proved that the Zariski tangent space to $R(G, \SL)$ at an irreducible representation $\rho$ can be identified to a subspace of $Z^1(G; \sll_{Ad \circ \rho})$. Here $\sll_{Ad \circ \rho}$ denotes the left $\ZZ[G]$-module structure on $\sll$ induced by the adjoint representation. This inclusion is explicitly given by \begin{equation}\label{InclusionNaturelle}
\TZar{\rho}{R(G, \SL)}  
\to Z^1(G; \sll_{Ad \circ \rho}),\: \diff{\rho_t}{t}{t=0} 	 
\mapsto \begin{cases}
   G \to \sll   & \text{ } \\
   g \mapsto \Diff{\rho_t(g)\rho(g^{-1})}{t}{t=0}   & \text{}
\end{cases}
\end{equation}
where $\rho_0 = \rho$ (\cf~\cite[Paragraph 3.1.3]{Porti:1997}). For an irreducible representation $\rho : G \to \SL$ the orbit $\SL(\rho) = \{Ad_A \circ \rho \; |\; A \in \PSL\}$ is isomorphic to  
$B^1_\rho(G; \sll_{Ad \circ \rho})$. So we get the natural inclusion
\[
T^{\mathrm{Zar}}_{\chi_\rho} X(G) \to Z^1(G; \sll_{Ad \circ \rho})/B^1(G; \sll_{Ad\circ \rho}) = H^1(G; \sll_{Ad \circ \rho}).
\]

	After this general digression we turn back to the case of $\SU$-representations. Consider the involution $\sigma : \SL \to \SL$ defined by $\sigma(A) = (\bar{A}^T)^{-1}$. Let $\SL^\sigma$ denote its set of fixed points. With this notation we have $\SL^\sigma = \SU$. Let $\sigma_*$ be the linear map $D_{\I} \sigma : \sll \to \sll$. We have $\sigma_*(x) = -\bar{x}^T$ and $\sll^{\sigma_*} = \su$. For every irreducible representation $\rho : G \to \SU$ we define the {tangent space} to $R(G) = (R(G, \SL))^{\sigma}$ at $\rho$ to be $$\tangent{\rho}{R(G)} = \left( \TZar{\rho}{R(G, \SL)} \right)^{\sigma_*} \subset \TZar{\rho}{R(G, \SL)}.$$

	The equality $Z^1(G; \sll_{Ad \circ \rho})  = Z^1_\rho(G) \otimes_{\IR} \CC$ provides the inclusion
\begin{equation}\label{TRintoDer}
\tangent{\rho}{R(G)} = \left( \TZar{\rho}{R(G, \SL)} \right)^{\sigma_*} \hookrightarrow  \left( Z^1(G; \sll_{Ad \circ \rho}) \right)^{\sigma_*} = Z^1_\rho(G).
\end{equation}

We also have: 
\begin{equation}\label{equality}
(\SL(\rho))^{\sigma} = \SU(\rho) \text{ for all } \rho \in \widetilde{R}(G).
\end{equation} 
Equality~(\ref{equality}) follows from the two following facts:
\begin{enumerate}
  \item two irreducible representations of $G$ in $\SL$ with the same character are conjugate by an element of $\SL$, see \cite[Proposition 1.5.2]{CS:1983};
  \item if two irrreducible (\ie non abelian) representations in $\SU$ are conjugate by an element of $\SL$, then they are conjugate by an element of $\SU$, see~\cite[Proof of Proposition 15]{Klassen:1991}.
\end{enumerate}

Inclusion~(\ref{TRintoDer}) and equality~(\ref{equality}) combine to yield 

\begin{prop}\label{isoRandH}
If $G$ is the free group $F_n$ (resp. the fundamental group of a punctured $2$-sphere), then inclusion~(\ref{TRintoDer}) gives rise to an isomorphism $\xymatrix@1@-.7pc{\tangent{[\rho]}{\widehat{R}(G)} \ar[r]^-\cong& H^1_\rho(G)}$ for every irreducible representation $\rho : G \to \SU$.
\end{prop}

\section{Twisted cohomology for knot exteriors}\label{TwistH}
	
	In this section, we turn to the geometric application of the previous algebraic preliminaries. 
	
	Assume that $S^3$ and $K$ are {oriented}. Let $\rho$ be any {non boundary central} representation of $G_K$ in $\SU$, \ie such that $\rho(\pi_1(\bord M_K)) \not\subset \{\pm \I\}$, and fix a presentation of the group of $K$ of the form:
\begin{equation}\label{PresW}
G_K = \langle S_1, \ldots, S_r \;|\; R_1, \ldots, R_{r-1} \rangle.
\end{equation}
The meridian $m$ of $K$ is only defined up to conjugation and is oriented by the convention $\lk(K,m) = +1$, where $\lk$ denotes the linking number. 

We begin by reviewing the notions of {regularity} and of {$\mu$-regularity} for a re\-pre\-sen\-ta\-tion (see~\cite{BZ:1996},~\cite{HeuKlassen:1997}, \cite[Definition 3.21]{Porti:1997} and~\cite[Section 1]{Heu:2003}).  

\subsection{Regular representations}
\label{RegularRep}
	The long exact sequence in $(Ad \circ \rho)$-twisted cohomology corresponding to the pair $(M_K, \bord M_K)$ and the Poincaré duality imply $\dim H^1_\rho(M_K) \geqslant 1$. So a $\SU$-representation of a knot group is {never} acyclic.
	
	If $\rho : G_K \to \SU$ is irreducible, then $H^0_\rho(M_K) =0$. Moreover, we have $\dim H^2_\rho(M_K) = \dim H^1_\rho(M_K)$ because $\sum_i (-1)^i \dim H^i_\rho(M_K) = 3\chi(M_K) = 0$. 

\subsubsection*{\textbf{Regular representations}}
	Among irreducible representations we focus on the regular ones. An irreducible $\SU$-representation $\rho$ of $G_K$  is called {regular} if $\dim H^1_\rho(M_K) = 1$. 
It is easy to see that this notion is invariant under conjugation. We let $\Reg{K}$ denote the set of conjugacy classes of regular representations of $G_K$ in $\SU$.

\begin{example}
	If $K$ denotes a torus knot or the figure eight knot, then any irreducible representation of $G_K$ in $\SU$ is regular (see the proof of Proposition~\ref{propEx}).
\end{example}
		
	M. Heusener and E. Klassen proved that if $\rho$ is regular, then its conjugacy class $[\rho]$ is a smooth point of $\widehat{R}(M_K)$ and $\dim \widehat{R}(M_K) = 1$ in a neighbourhood of $[\rho]$ (\cf~\cite[Proposition 1]{HeuKlassen:1997}). So, $\Reg{K}$ is a {$1$-dimensional manifold}. Using properties of Subsection~\ref{DerTan}, we can establish 

\begin{prop}\label{IsomTH}
If $\rho : G_K \to \SU$ is regular, then inclusion~(\ref{TRintoDer}) gives rise to an isomorphism $\varphi_{[\rho]} : \xymatrix@1@-.7pc{\tangent{[\rho]}{\widehat{R}(M_K)} \ar[r]^-\cong &H^1_\rho(M_K)}$.
\end{prop}

\subsubsection*{\textbf{$\mu$-regular representations}}
	Consider $\rho \in \widetilde{R}(M_K)$ and fix a simple closed oriented curve $\mu$ in $\bord M_K$.  Let $N_{[\rho]}$ be a neighbourhood of $[\rho]$ in $\widehat{R}(M_K)$. Introduce the map $\theta_\mu : N_{[\rho]}\subset \widehat{R}(M_K) \to S^1$ which associates to $[\varrho] \in N_{[\rho]}$ half the angle of the rotation $Ad \circ \varrho(\mu)$. As we only consider knot exteriors remark that the map $\theta_\mu$ is never zero and always a well-defined analytic map.
	
	A representation $\rho$ is called {$\mu$-regular} if $\rho$ is regular and if $\theta_\mu : \widehat{R}(M_K) \to S^1$ is a submersion at $[\rho]$ (see~\cite[Proposition 3.26]{Porti:1997}). Observe that the notion of $\mu$-regularity is invariant under conjugation and does not depend on the orientation of $\mu$.

\begin{example}
	If $K$ denotes a torus knot, then any irreducible representation of $G_K$ in $\SU$ is $m$-regular (see~\cite[Example 1.43]{JDTHESE} and the proof of Proposition~\ref{propEx}).
\end{example}

	We use in the sequel the following alternative cohomology formulation of $\mu$-regu\-la\-ri\-ty. Consider the linear form $f^\rho_\mu : H^1_\rho(M_K) \to \IR$ defined by $f^\rho_\mu(v) = \langle v(\mu), P^\rho \rangle$. Here $\langle \cdot, \cdot \rangle$ denotes the usual scalar product of $\su$.
	
\begin{prop}\label{propmureg}
The representation $\rho$ is $\mu$-regular if and only if $f^\rho_\mu$ is an isomorphism. 
\end{prop}
It follows that if $\rho$ is $\mu$-regular, then we have a reference generator for $H^1_\rho(M_K) \cong \IR$ denoted $h^{(1)}_\rho(\mu)$ which satisfies $f^\rho_\mu(h^{(1)}_\rho(\mu)) = +1$.

Proposition~\ref{propmureg} is a consequence of the following claim.
\begin{claim}\label{claimreg}
	If $\rho \in R(M_K)$ is regular, then $f^\rho_\mu(v) = D_{[\rho]}\theta_\mu(\varphi^{-1}_{[\rho]}(v))$, for all $v \in H^1_\rho(M_K)$. 
\end{claim}
\begin{proof}
	Fix $v \in H^1_\rho(M_K)$ and choose a germ $\rho_s$ in a neighbourhood of the origin such that $\varphi^{-1}_{[\rho]}(v) = \diff{\rho_s}{s}{s=0}$ with $\rho = \rho_0$. We can assume that $\rho_s$ is regular because $\rho$ is itself regular. There exists a family of matrices $A_s$ such that
\[
\rho_s(\mu) = A_s  \left(\begin{array}{cc}e^{i \theta_\mu(\rho_s)} & 0 \\0 & e^{-i \theta_\mu(\rho_s)}\end{array}\right) A_s^{-1}.
\]
Observe that $P^\rho = Ad_{A_0}(\ii)$ and $\diff{A_s^{-1}}{s}{s=0} = -A_0^{-1} \diff{A_s}{s}{s=0} A_0^{-1}.$ Thus, 
\begin{equation}\label{diffrho}
\Diff{\rho_s(g)\rho(g)^{-1}}{s}{s=0} = (\mathrm{Id} - Ad_{\rho(\mu)}) \diff{A_s}{s}{s=0} A_0^{-1}+ \diff{\theta_\mu({\rho_s})}{s}{s=0} P^\rho.
\end{equation}
Taking into account the fact that $(\mathrm{Id} - Ad_{\rho(\mu)}) P^\rho = 0$ we deduce that the term $(\mathrm{Id} - Ad_{\rho(\mu)}) \diff{A_s}{s}{s=0} A_0^{-1}$ in equation~(\ref{diffrho}) is orthogonal to $P^\rho$. Thus,
\[
f^\rho_\mu(v) = \left\langle \Diff{\rho_s(g)\rho(g)^{-1}}{s}{s=0}, P^\rho \right\rangle = \diff{\theta_\mu({\rho_s})}{s}{s=0}.
\]
\end{proof}

\subsection{Cohomology orientation of knot exteriors}
\label{CohomOr}
In this subsection, we equip knot exteriors with the {canonical cohomology orientation} (see~\cite[Chapter V.3]{Turaev:2002} for more details).

	We associate to the presentation~(\ref{PresW}) of $G_K$ the two-dimensional CW-com\-plex $X_K$ in the usual way. The 0-skeleton of $X_K$ consists of a single point $e^{(0)}$, the $1$-skeleton $X_K^1$ is a wedge of $r$ oriented circles $e^{(1)}_1, \ldots, e^{(1)}_r$ corresponding to the generators $S_i$. Finally, $X_K$ is obtained from $X_K^1$ by gluing the $r-1$ closed $2$-cells $e^{(2)}_1, \ldots, e^{(2)}_{r-1}$ attached using the relations $R_j$. We have $\pi_1(X_K) = G_K$. In~\cite{Waldhausen:1978}, F. Waldhausen proved that the {Whitehead group} of a knot group is {trivial}. As a result $X_K$ has {the same simple homotopy type} as $M_K$. So, the CW-complex $X_K$ can be used to compute the Reidemeister torsion of $M_K$. 
	
	Set $\mathfrak{X}_i = C^{2-i}(X_K; \IR)$ for  $i=0,1,2$. For consistency reasons, we think of $\mathfrak{X}_*$ as a chain complex of real vector spaces. Explicitly
\begin{equation}\label{NoeudReel}
\mathfrak{X}_* = \xymatrix@1{0 \ar[r] & \IR e^{(0)}_{} \ar[r] & \displaystyle{\bigoplus_{j=1}^{r} \IR e^{(1)}_j} \ar[r] & \displaystyle{\bigoplus_{k=1}^{r-1}\IR e^{(2)}_k} \ar[r] & 0}.
\end{equation}
Observe that $\chi(\mathfrak{X}_*) = 0$ and that $\mathfrak{X}_*$ is not acyclic: $H_i(\mathfrak{X}_*) \cong H^{2-i}(M_K; \IR)$, for $i=0, 1, 2$.

	Let $\lbrack \! \lbrack pt \rbrack \! \rbrack \in H^0(M_K; \ZZ)$ be the cohomology class of a point and let $m^* : m \mapsto 1$ be the dual of $m$. We base $H_*(\mathfrak{X}_*) = H^0(M_K; \IR) \oplus H^1(M_K; \IR)$ with $\left\{ \lbrack \! \lbrack pt \rbrack \! \rbrack, m^*\right\}$. In the sequel, we  assume that $\mathfrak{X}_*$ endows the cohomology orientation induced by this basis. Moreover we always compute the $(Ad \circ \rho)$-twisted Reidemeister torsion  of knot exteriors with respect to this canonical cohomology orientation.

\subsection{The knot exterior twisted chain complex}
\label{TwistedC}
	Let $\rho$ be a non boundary central representation of $G_K$. The \emph{knot exterior twisted chain complex} is the chain complex $\mathscr{X}^\rho_*$ defined by $\mathscr{X}^\rho_i = C^{2-i}(X_K; Ad \circ \rho)$ for $i=0,1,2$. As in Subsection~\ref{CohomOr}  we think of $\mathscr{X}^\rho_*$ as a chain complex of real vector spaces. The presentation~(\ref{PresW}) of $G_K$ provides 
\begin{equation}\label{twistX}
\mathscr{X}^\rho_* =  \xymatrix@1{
0 \ar[r] & \su \ar[r]^-{d^{\rho}_2} & \su^{r} \ar[r]^-{d^{\rho}_1} & \su^{r-1} \ar[r] & 0}.
\end{equation}
The boundary operators $d^\rho_1$ and $d^\rho_2$ are induced by the usual ones of $C_*(\widetilde{M}_{K}, \ZZ)$ obtained using the Fox differential calculus (see~\cite[Chapter 7]{Crowell-Fox}).  If we write  $g \circ x = Ad_{\rho(g)}(x)$, for $g \in G_K$ and $x \in \su$, then 
$$
d^\rho_2(x) = \left( (1-S_1) \circ x, \ldots, (1-S_r) \circ x\right) \text{ for all } x \in \su,
$$
\[
d^\rho_1 \left( {(x_j)_{1 \leqslant j\leqslant r}} \right) = \left( \sum_{j=1}^r \fox{R_i}{S_j} \circ x_{j}\right)_{1\leqslant i \leqslant r-1} \text{ for all } (x_j)_{1 \leqslant j\leqslant r} \in \su^{r}.
\]

\subsection{The Reidemeister torsion and the Alexander polynomial}
\label{SS:polynomeAlexander}
	In the previous subsections, we focus on non abelian representations. In the present one, let us make a digression on abelian representations. 
	The aim of this subsection is to compute the Reidemeister torsion of the exterior of $K$ twisted by the adjoint representation associated to an abelian representation of $G_K$ in terms of the Alexander polynomial of $K$.
	
	Recall that an abelian representation of $G_K$ in $\SU$ is entirely determined by its value on the meridian $m$ of $K$. 
	Let $\varphi_\theta : G_K \to \SU$ be the abelian representation defined by $\varphi_\theta(m) = \cos(\theta) + \sin(\theta)\ii$ with $0 < \theta < \pi$ (\ie $\varphi_\theta$ is supposed to be non boundary central). When $e^{2i\theta}$ is {not} a zero of the Alexander polynomial $\Delta_K$ of $K$ we say that $\varphi_\theta$ is {regular}. In this case, E. Klassen proved in~\cite[Theorem 19]{Klassen:1991} that $H^i_{\varphi_\theta}(M_K) \cong H^i(M_K; \ZZ) \otimes \IR$, for all $i$. 
	We can also remark that $H^0_{\varphi_\theta}(M_K)=\ker d_2^{\varphi_\theta}$ is generated by $h^{(0)} = \ii$ (the common fixed axis of the rotations $Ad \circ \varphi_{\theta}(S_j)$, $1 \leqslant j \leqslant r$) and $H^1_{\varphi_\theta}(M_K)$ is generated by $h^{(1)} = \ii_1 + \cdots + \ii_{2n}$, where $\ii_k$ is the vector in $\su^{2n}$ of which all entries are zero except the one of index $k$ which is equal to $\ii$.
	
	With this notation we have (see~\cite[Theorem 4]{Milnor:1962}, \cite[Subsection 1.2]{Turaev:1986} and~\cite[Section 2.5]{JDTHESE})
	
\begin{prop}\label{P:TorsionAlexander}
	Let $\varphi_\theta$ be a regular abelian representation with $0<\theta<\pi$. The $(Ad \circ \varphi_\theta)$-twisted Reidemeister torsion of $M_K$ calculated in the basis $\{ h^{(0)}, h^{(1)}\}$ of $H^*_{\varphi_\theta}(M_K)$ and with respect to the canonical cohomology orientation of $M_K$ satisfies
\begin{equation}
\label{EQ:polynomeAlexander}
\mathrm{TOR}\left({M_K; Ad \circ \varphi_\theta, \{h^{(0)}, h^{(1)}\}} \right) = \frac{4\sin^2(\theta)}{|\Delta_K(e^{2i \theta})|^2}.
\end{equation}
Here $\Delta_K$ denotes the Alexander polynomial of $K$.
\end{prop}
\begin{remark} If $K$ is the trivial knot, then $$\mathrm{TOR}\left({M_K; Ad \circ \varphi_\theta, \{h^{(0)}, h^{(1)}\}}\right) = 4\sin^2(\theta)$$ is the twisted Reidemeister torsion of the solid torus $M_K$.
\end{remark}
\begin{proof}[Ideas of the proof] We just give the main ideas and leave the detailed computations. For all $j$, in the basis $\{\ii, \jj, \kk\}$, we have $$Ad\circ \varphi_\theta (S_j) = \begin{pmatrix}
	1 & 0 & 0 \\	0 & \cos(2\theta) & -\sin(2\theta) \\
	0 & \sin(2\theta) & \cos(2\theta) \\
\end{pmatrix} \in \mathrm{Aut}(\su).$$ 
Take the complexification of the Lie algebra $\su$. The three subspaces res\-pec\-ti\-ve\-ly generated by $\ii$, $\jj - i \kk$ and $\jj + i \kk$ define a splitting of $\su$. The action of $Ad\circ \varphi_\theta$ leaves these subspaces invariant: it acts trivially on the first factor, as the multiplication by $e^{2i\theta}$ on the second one and as the multiplication by $e^{-2i\theta}$ on the third one. The factor corresponding to the subspace generated by $\ii$ con\-tri\-bu\-tes trivially because $H_1(M_K)$ is torsion free. To compute the other contributions recall that the Rei\-de\-meister torsion twisted by the abelianization representation is (with our conventions for the Reidemeister torsion) equal to $({t-1})/{\Delta_K(t)}$, see~\cite{Milnor:1962}. This result implies that the factor corresponding to the subspace generated by $\jj - i \kk$ contributes by $({e^{2i\theta} - 1})/{\Delta_K(e^{2i\theta})}$ and that the factor corresponding to the subspace generated by $\jj + i \kk$ contributes by $({e^{-2i\theta} - 1})/{\Delta_K(e^{-2i\theta})}$. The product of these three contributions gives equation~(\ref{EQ:polynomeAlexander}).
\end{proof}

\section{Non abelian Reidemeister torsion}\label{Torsion}

	In this section, we investigate the non abelian Reidemeister torsion for knots in $S^3$. 
	For this purpose, we define a  reference basis for the second twisted cohomology group $H^2_\rho(M_K)$. The main idea of the construction is to look at the restriction of the $\SU$-representations of $G_K$ to the peripheral subgroup of $K$ (see Subsections~\ref{TwistT} \&~\ref{RefBasis2}). 
	We assume that $S^3$ and $K \subset S^3$ are {oriented}. Later we shall see that the definition of the torsion does not depend on the orientation of $K$ but depends on the one of $S^3$ (see Proposition~\ref{P:dependanceTorsion}). 
	
\subsection{Twisted cohomology of the torus}\label{TwistT}
 
 	The exterior of $K$ is oriented and is bounded by the $2$-torus $\bord M_K$. This boundary inherits  an orientation by the convention ``\emph{the inward pointing normal vector in the last position}". 
	Let $\mathrm{int}(\cdot, \cdot)$ be the intersection form on $\bord M_K$ induced by this orientation. 
	The {peripheral subgroup} $\pi_1(\bord M_K)$ is generated by the meridian-longitude system $m, l$ of $K$. Here $m$ is oriented by the convention $\lk(K, m) = +1$ and $l$ is oriented by using the rule $\mathrm{int}(m, l) = +1$. 

	Observe that if $\rho \in \widetilde{R}(G_{K})$ then $\rho(m) \ne \pm \I$ (because $m$ normally generates $G_{K}$). So there exists a unique pair $(\theta, P^\rho) \in (0, \pi) \times S^2$ such that $\rho(m)=\cos(\theta) + \sin(\theta) P^\rho$. Notice that $P^\rho$ generates the common fixed axis of the rotations $Ad \circ \rho(\pi_1(\bord M_K))$ and thus generates the $1$-dimensional space $H^0_\rho(\bord M_K)$.

	 The usual scalar product of $\su$ induces a {cup-product}
$$\cup  \;: H^p_\rho(\bord M_K) \times H^{q}_\rho(\bord M_K) \to H^{p+q}(\bord M_K; \IR).$$ This cup-product is non-degenerated and gives an explicit isomorphism between $H^*_\rho(\bord M_K)$ and $H^*(\bord M_K; \IR)$ (see~\cite[Proposition 3.18]{Porti:1997} for a proof). We summarize this in the following result: 
\begin{lemma}\label{L:CohomologieTore}
The map $\phi^*_{P^\rho} : H^*_\rho(\bord M_K) \to H^*(\bord M_K; \IR)$, defined by $\phi^*_{P^\rho}(z) = P^\rho \cup z$, is a natural isomorphism.
\end{lemma}

\subsection{Reference basis of the second twisted cohomology group of $M_K$}
\label{RefBasis2}
	Recall that $H^1_\rho(M_K) \cong \IR$ and $H^2_\rho(M_K) \cong \IR$ for all regular representation $\rho$ of $G_K$ in $\SU$. The construction of the reference generator of $H^2_\rho(M_K)$ is based on the following lemma (see~\cite[Corollary 3.23]{Porti:1997}).
\begin{lemma}\label{L:IsomH2}
	Let $\rho : G_K \to \SU$ be a regular representation. The inclusion $\bord M_K \hookrightarrow M_K$ induces a natural isomorphism $i^*: H^2_\rho(M_K) \to H^2_\rho(\bord M_K)$.
\end{lemma}
\begin{proof}
	Associated to the pair $(M_K, \bord M_K)$ is the long exact sequence in $(Ad\circ \rho)$-twisted cohomology
\[
\xymatrix@1@-.6pc{\cdots \ar[r] & H^1_\rho(\bord M_K) \ar[r] & H^2_\rho(M_K, \bord M_K) \ar[r] & H^2_\rho(M_K) \ar[r]^-{i^*} & H^2_\rho(\bord M_K) \ar[r] & 0.}
\]
Thus 
$\mathrm{rk}\, i^* = \dim H^2_\rho(\bord M_K) = \dim H^2(\bord M_K; \IR) = 1$ (see Lemma~\ref{L:CohomologieTore}). Furthermore $\dim H^2_\rho(M_K) = 1$. So $i^* $ is an isomorphism.
\end{proof}

	Joining together Lemmas~\ref{L:CohomologieTore} \&~\ref{L:IsomH2} we conclude that the composition $$\phi^{(2)}_{P^\rho} \circ i^* : H^2_\rho(M_K) \to H^2_\rho(\bord M_K) \to H^2(\bord M_K; \IR)$$ is an {isomorphism}. Let $c$ denote the generator of $H^2(\bord M_K; \ZZ)$ corresponding to the fundamental class $\lbrack \! \lbrack \bord M_K \rbrack \! \rbrack \in H_2(\bord M_K; \ZZ)$ induced by the orientation of $\bord M_K$. The \emph{reference generator} $h^{(2)}_\rho$ of $H^2_\rho(M_K)$ is defined by 
\begin{equation}\label{EQ:Defh2}
h^{(2)}_\rho = (\phi^{(2)}_{P^\rho} \circ i^*)^{-1}(c).
\end{equation}

\subsection{Definition of the volume form $\tau^K$}
\label{defTau}

	Let $\rho \in R(G_K)$ be regular. The \emph{Rei\-de\-meis\-ter torsion form $\tau^K$} at $\rho$ is the linear form $\tau^{K}_\rho: \tangent{[\rho]}{\widehat{R}(M_K)} \to \mathbb{R}$  
defined by $$\tau^{K}_\rho(v) = \mathrm{TOR}\left(M_K; Ad \circ \rho, \left\{{\varphi_{[\rho]}(v), h^{(2)}_\rho}\right\} \right),  \text{ for all } v \ne 0.$$ Here $\varphi_{[\rho]} : \xymatrix@1@-.7pc{\tangent{[\rho]}{\widehat{R}(M_K)} \ar[r]^-{\cong} &H^1_\rho(M_K)}$  (see Proposition~\ref{IsomTH}) and $h^{(2)}_\rho$ is the reference generator of  $H^2_\rho(M_K)$ (see equation~(\ref{EQ:Defh2})).

	With the notation of Subsections~\ref{CohomOr}--\ref{TwistedC} we observe $$\tau^K_\rho(v) = \mathrm{sgn}(\mathrm{Tor}(\mathfrak{X}_*)) \cdot \mathrm{Tor}\left(\mathscr{X}^\rho_*, \left\{{\varphi_{[\rho]}(v), h^{(2)}_\rho}\right\}\right).$$
	
	The linear form $\tau^K_\rho$ just depends on the $\SO$-conjugacy class $[\rho]$ of $\rho$, thus $\tau^K : [\rho] \mapsto \tau^K_\rho$ is a well-defined $1$-volume form on $\Reg{K}$.
	
Here are some remarks.
\begin{remark}\label{remsign}
Note that $\mathrm{Tor}(\mathfrak{X}_*)$ is $\pm 1$, because the homology of $M_K$ with integral coefficients is torsion free (see~\cite{Turaev:1986}). So that $\mathrm{sgn}(\mathrm{Tor}(\mathfrak{X}_*)) = \mathrm{Tor}(\mathfrak{X}_*)$.
\end{remark}
	
\begin{remark}\label{Rem}
	For a regular representation $\rho$, there exist a unique $P^\rho \in S^2$ and a unique $\bar{P}^\rho \in S^2$ such that:
$\rho(m) = \cos(\theta) + \sin(\theta) P^\rho$  and $\rho(m) = \cos(2\pi-\theta) + \sin(2\pi-\theta) \bar{P}^\rho$  where $0 < \theta < \pi.$
If we repeat the construction of the torsion form replacing the vector $P^\rho$ by  $\bar{P}^\rho$, then we obtain the opposite volume form $- {\tau}^K_\rho$ (because $h^{(2)}_\rho$ is replaced by $-h^{(2)}_\rho$, see Lemma~\ref{L:CohomologieTore}). 
\end{remark}

\begin{prop}\label{P:dependanceTorsion}
	If $K \subset S^3$ is a knot, then ${\tau}^K$ has the following properties:
\begin{enumerate}
  \item ${\tau}^K$ does not depend on the orientation of $K$.
  \item If $K^*$ denotes the mirror image of $K$, then we have $$(\Reg{K^*}, {\tau}^{K^*}) = (-\Reg{K}, -{\tau}^{K}).$$
\end{enumerate}
\end{prop}
\begin{proof} It requires two steps.
\begin{enumerate}
  \item If we reverse the orientation of $K$, then the meridian $m$ is changed in $m^{-1}$ and the orientation of $M_{K}$ is not affected. So that the cohomology orientation is reversed and the reference generator $h^{(2)}_\rho$ defined by~(\ref{EQ:Defh2}) is replaced by $-h^{(2)}_\rho$. Formula~(\ref{EQ:changementdebase}) implies that $\mathrm{Tor}(\mathfrak{X}_*)$ and $\mathrm{Tor}(\mathscr{X}^{\rho}_*)$ change their sign at the same time. 
  \item If we reverse the orientation of $S^{3}$, then the orientation of $M_{K}$ and the one of $\bord M_{K}$ change at the same time. Hence the meridian $m$ is changed in $m^{-1}$. So that the cohomology orientation is reversed. Thus $\mathrm{Tor}(\mathfrak{X}_*)$ changes its sign. Moreover  the reference generator $h^{(2)}_\rho$ defined by~(\ref{EQ:Defh2}) is unchanged. Thus $\mathrm{Tor}(\mathscr{X}^{\rho}_*)$ does not change. Moreover 
\end{enumerate}
\end{proof}

\subsection{An explicit computation}

	Let $q>0$ be an odd integer. Let $K_q$ denote the torus knot of type $(2, q)$ and let $M_q$ denote its exterior. Recall that the group $G_q$ of $K_q$ admits the presentation $G_q = \langle x, y \; |\; x^2 = y^q\rangle$. In $G_q$ the meridian $m$ of $K_q$ can be written as $m = xy^{\frac{1-q}{2}}$ and the longitude $l$ as $l =  x^2m^{-2q}$. 
	
	Before making the explicit computation of $\tau^{K_q}$ we give a parametrization of $\widehat{R}(M_q)$ (cf.~\cite[Theorem 1]{Klassen:1991}). Consider the map $\rho_{\ell,t} : G_q \to \SU$ defined by $\rho_{\ell,t}(x) = \ii$ and $\rho_{\ell,t}(y)=\cos\left( {(2\ell-1)\pi}/{q}\right)+\sin\left( {(2\ell-1)\pi}/{q}\right)\left( \cos(\pi t)\ii+\sin(\pi t)\jj \right),$ for $0<t<1$ and $\ell \in \left\{1, \ldots, {(q-1)}/{2}\right\}.$
	For all $t$ and $\ell$, $\rho_{\ell, t}$ defines an irreducible representation of $G_q$. In fact, {all} irreducible representations of $G_q$ in $\SU$ arise in this way; more precisely, any $\rho \in \widetilde{R}(M_q)$ is conjugate to one and only one of the $\rho_{\ell,t}$. With this notation we show
	
\begin{prop}\label{propEx}
	We have $\Reg{K_q} = \widehat{R}(M_q)$. The Reidemeister torsion form associated to the right-hand torus knot $K_q$ of type $(2, q)$ satisfies
\[
\tau^{K_q}_{[\rho_{\ell, t}]}\left(\diff{\rho_{\ell,s}}{s}{s=t}\right) = - \frac{8}{q} \sin^2 \left( \frac{(2\ell-1)\pi}{q} \right) \diff{\theta_m^{\rho_{\ell, s}}}{s}{s=t}
\] 
where $\theta_m^{\rho_{\ell, t}}= \theta_m({\rho_{\ell, t}}) = \arccos\left((-1)^{\ell-1}\cos\left((2\ell-1)\pi/2q\right)\cos(\pi t) \right).$
\end{prop}
\begin{proof}
	The proof is a direct computation. Essentially it consists in the explicit determination of the generator $h^{(2)}_{\rho_{\ell, t}}$ of $H^2_{\rho_{\ell, t}}(M_q)$ defined by equation~(\ref{EQ:Defh2}).
	
\noindent (1) For all $t$ and $\ell$, observe that $\rho_{\ell, t}(m) = \cos(\theta_m^{\rho_{\ell, t}}) + \sin(\theta_m^{\rho_{\ell, t}}) P^{\rho_{\ell, t}}$ with 
\[
P^{\rho_{\ell, t}} = (-1)^{\ell-1} \frac{\sin\left( \frac{(2\ell-1)\pi}{2q}\right)\ii - \cos\left( \frac{(2\ell-1)\pi}{2q}\right)\sin(\pi t)\kk}{\sqrt{\sin^2\left( \frac{(2\ell-1)\pi}{2q}\right)+\cos^2\left( \frac{(2\ell-1)\pi}{2q}\right)\sin^2(\pi t)}}.
\]

\noindent (2) Let $X_q$ denote the $2$-dimensional CW-complex associated to the presentation $G_q = \langle x, y \;|\; r\rangle$, where $r=x^2y^{-q}$.  The corresponding twisted chain complex $\mathscr{Z}^{\rho_{\ell, t}}_* = C^{2 - *}(X_q; Ad \circ {\rho_{\ell, t}})$ is 
\[
\mathscr{Z}^{\rho_{\ell, t}}_* = \xymatrix{0 \ar[r] & \su \ar[r]^-{d_2^{\rho_{\ell, t}}} & \su \oplus \su \ar[r]^-{d_1^{\rho_{\ell, t}}} & \su \ar[r] & 0}
\]
and the coboundary operators are 
\[
d_2^{\rho_{\ell, t}} = \binom{Ad_X-\mathrm{Id}}{Ad_Y-\mathrm{Id}}  \text{ and }  d_1^{\rho_{\ell, t}} = \left(\mathrm{Id}+Ad_X, -(\mathrm{Id}+Ad_Y+\cdots+Ad_{Y^{q-1}})\right).
\]
Here $X =\rho_{\ell, t}(x)$, $Y = \rho_{\ell, t}(y)$. Observe that $\coker d_1^{\rho_{\ell, t}} \cong \IR$ is generated by $\kk$. So $\dim H^2_{\rho_{\ell, t}}(M_q) = 1$ and thus $\rho_{\ell, t}$ is regular.

	Now, we compute the homomorphism $i^* : H^2_{\rho_{\ell, t}}(M_q) \to H^2_{\rho_{\ell, t}}(\bord M_q)$ induced by $i : \bord M_q \hookrightarrow M_q$. For this purpose, we write the relation of $\pi_1(\bord M_q)$ as follows: $$R = lml^{-1}m^{-1} = xrx^{-1}(xy^{\frac{1-q}{2}})r^{-1}(xy^{\frac{1-q}{2}})^{-1}.$$ So that the inclusion $i$ induces, at the level of $2$-cocycles, the map $i^* : \su \to \su$ defined by $i^*(z) = Ad_X(z) - Ad_{W}(z)$, where $W = \rho_{\ell, t}(m)$. Thus, the re\-fe\-ren\-ce generator $h^{(2)}_{\rho_{\ell, t}}  \in H^2_{\rho_{\ell, t}}(M_q)$ is the equivalence class of the element $u^{(2)} = {\langle Ad_X(\kk) - Ad_{W}(\kk), P^{\rho_{\ell, t}}\rangle}^{-1} \kk \in \su$.
	
\noindent (3) The computation of the torsion for the torus knot $K_q$ requires a reference generator for $H^1_{\rho_{\ell, t}}(M_q)$.
	Observe that $v=\binom{\kk}{0} \in \ker d_1^{\rho_{\ell, t}} \setminus \im  d_2^{\rho_{\ell, t}}$. Furthermore $v(m) = v(x) = \kk$, so that the equivalence class of the element $u^{(1)} = {\langle \kk, P^{\rho_{\ell, t}} \rangle}^{-1} v$ generates $H^1_{\rho_{\ell, t}}(M_q)$.

\noindent (4)  A direct determinant computation gives $$\mathrm{Tor}(\mathscr{Z}_*^{\rho_{\ell, t}}, \{u^{(1)}, u^{(2)}\}) = -\frac{8}{q}\sin^2\left( \frac{(2\ell-1)\pi}{q}\right).$$

	For the right-hand torus knot we have $\tau_0 = \mathrm{sign}(\mathrm{Tor}(C^{*}(X_q; \IR))) = +1.$  Thus 
$$\mathrm{TOR}(M_q ; Ad \circ \rho_{\ell, t}, \{u^{(1)}, u^{(2)}\}) = \tau_0 \mathrm{Tor}(\mathscr{Z}_*^{\rho_{\ell, t}}, \{u^{(1)}, u^{(2)}\}) = -\frac{8}{q}\sin^2\left( \frac{(2\ell-1)\pi}{q}\right).$$
	
	The basis change formula~(\ref{EQ:changementdebase}) provides $$\tau^{K_q}_{[\rho_{\ell, t}]}\left(\diff{\rho_{\ell,s}}{s}{s=t}\right) = \left[\diff{\rho_{\ell,s}}{s}{s=t} / u^{(1)} \right] \cdot \mathrm{TOR}(M_q ; Ad \circ \rho_{\ell, t}, \{u^{(1)}, u^{(2)}\}),$$
and because
\[
\left[\diff{\rho_{\ell,s}}{s}{s=t} / u^{(1)} \right] =\left\langle \Diff{\rho_{\ell,s}(m)\rho_{\ell,t}(m)^{-1}}{s}{s=t} , 								P^{\rho_{\ell, t}}\right\rangle
=\diff{\theta_m^{\rho_{\ell, s}}}{s}{s=t}
\]
we obtain the required formula of the proposition.
\end{proof}

\begin{remark}
	Proposition~\ref{propEx} can be deduced from the computations made by J. Park for graph manifolds in~\cite{Park:1997}. However the distinguished generator of $H^2_{\rho_{\ell, t}}(M_q)$ must be explicitly known which is the main difficult part of the preceding proof.
\end{remark}

\section{A volume form associated to a knot}\label{Vol}
	 In this section, we review in detail the ``natural" Casson-type construction of the {volume form} on $\Reg{K}$, see \cite{CRAS} and~\cite{VolumeForm}.
	
\subsection{The ``base $\wedge$ fiber" condition} 
	A {volume form} $v$ on a $n$-dimensional manifold is a nowhere vanishing differential $n$-form. In the sequel, we will make a great use of the  \emph{``base $\wedge$ fiber" condition} which is the following. Given two volume forms $v$ and $w$ on the manifolds $M^m$ and $N^n$ respectively, a submersion $f : M \to N$ and a point $y \in N$, then the subspace $f^{-1}(y) \subset M$ is a  submanifold of dimension $m-n$, the tangent space $\tangent{x}{f^{-1}(y)}$ is the kernel of $D_x f$ and we have the short exact sequence
\begin{equation*}\label{Base+Fibre}
\xymatrix@1@-.5pc{0 \ar[r] & \tangent{x}{f^{-1}(y)} \ar[r]^-i & \tangent{x}{M} \ar[r]^-{D_xf} & \tangent{y}{N} \ar[r] & 0.}
\end{equation*}
The submanifold  $f^{-1}(y)$ is endowed with the unique volume form $\omega$ such that, for all $x \in f^{-1}(y)$, $\omega_x = v_x/w_y$ \ie $\omega_x \wedge w_y = (i \oplus s)^*(v_x)$, $s$ being a section of $D_xf$. 

\subsection{Plat decomposition and splitting of knot exteriors}
\begin{figure}[!htp]
\begin{center}
\begin{pspicture}(3,5.5)
%
\psline(1.5,3.5)(1.5,3.75)\psline(2,3.5)(2,3.75)
\psline(1.5,1.25)(1.5,1.5)\psline(2,1.25)(2,1.5)
\psline(1.5,3.5)(2,3)\psline(2,3.5)(1.8,3.3)\psline(1.7,3.2)(1.5,3)
\psline(1.5,3)(1.5,2.75)\psline(2,3)(2,2.75)
\psline(1.5,2.75)(2,2.25)\psline(2,2.75)(1.8,2.55)\psline(1.7,2.45)(1.5,2.25)
\psline(1.5,2.25)(1.5,2)\psline(2,2.25)(2,2)
\psline(1.5,2)(2,1.5)\psline(2,2)(1.8,1.8)\psline(1.7,1.7)(1.5,1.5)
\pscurve(.5,3.95)(.8,5.2)(1,5.3)(1.2,5.2)(1.5,3.95)
\pscurve(2,3.95)(2.3,5.2)(2.5,5.3)(2.7,5.2)(3,3.95)
%
\pscurve{<-}(.5,3.75)(.65,.45)(.8,.1)(1,0)(1.2,.1)(1.35,.45)(1.4,.7)(1.5,1.25)
\pscurve(2,1.25)(2.3,.1)(2.5,0)(2.7,.1)(3,1.25)(3,3.75)
%
%
\psline[linewidth=.005](4,3.85)(4.5,5)
%
\psline[linewidth=.005](-.8,3.85)(0,5)
\psline[linewidth=.005](-.8,3.85)(4,3.85)
\psline[linewidth=.005](0,5)(.6,5)
\psline[linewidth=.005](.8,5)(1.2,5)
\psline[linewidth=.005](1.45,5)(2.1,5)
\psline[linewidth=.005](2.3,5)(2.7,5)
\psline[linewidth=.005](2.9,5)(4.5,5)
%
%
\uput{0}[0](1.55,-.2){$B_1$}
\uput{0}[0](1.55,5.5){$B_2$}
\uput{0}[0](4,4.75){$S$}
%
\psline[linewidth=.1, linecolor=white](.5,4.075)(.5,4.2)
\psline[linewidth=.1, linecolor=white](1.5,4.075)(1.5,4.2)
\psline[linewidth=.1, linecolor=white](2,4.075)(2,4.2)
\psline[linewidth=.1, linecolor=white](3,4.075)(3,4.2)
\psline[linewidth=.1, linecolor=white](.775,5.15)(.85,5.3)
\psline[linewidth=.1, linecolor=white](2.65,5.3)(2.725,5.15)
\psline[linewidth=.1, linecolor=white](.95,0)(1.075,0)
\psline[linewidth=.1, linecolor=white](2.55,0)(2.425,0)
\psarc[linewidth=.005, linecolor=black]{->}(.9,5.255){.125}{65}{325}
\psline[linewidth=.005, linecolor=black](.85,5.355)(.5,5.75)
\psarc[linewidth=.005, linecolor=black]{->}(2.6,5.26){.125}{195}{105}
\psline[linewidth=.005, linecolor=black](2.675,5.35)(3,5.75)
\psarc[linewidth=.005, linecolor=black]{->}(.9,.04){.125}{195}{110}
\psline[linewidth=.005, linecolor=black](.85,-.06)(.5,-.35)
\psarc[linewidth=.005, linecolor=black]{->}(2.6,.04){.125}{65}{0}
\psline[linewidth=.005, linecolor=black](2.65,-.06)(3,-.35)
\psarc[linewidth=.005, linecolor=black]{->}(.5,4.25){.125}{115}{40}
\psline[linewidth=.005, linecolor=black](.45,4.15)(.3,3.85)
\psarc[linewidth=.005, linecolor=black]{->}(1.45,4.25){.125}{115}{40}
\psline[linewidth=.005, linecolor=black](1.4,4.15)(1.25,3.85)
\psarc[linewidth=.005, linecolor=black]{->}(2,4.25){.125}{115}{40}
\psline[linewidth=.005, linecolor=black](1.95,4.15)(1.8,3.85)
\psarc[linewidth=.005, linecolor=black]{->}(2.95,4.25){.125}{115}{40}
\psline[linewidth=.005, linecolor=black](2.9,4.15)(2.75,3.85)
%
%
\uput{0}[0](.3,5.35){\Tiny {$t_1^{(2)}$}}
\uput{0}[0](2.8,5.35){\Tiny {$t_2^{(2)}$}}
\uput{0}[0](.15,4.2){\Tiny {$s_1$}}
\uput{0}[0](1.05,4.2){\Tiny {$s_2$}}
\uput{0}[0](2.15,4.2){\Tiny {$s_3$}}
\uput{0}[0](3.1,4.2){\Tiny {$s_4$}}
\uput{0}[0](2.8,0){\Tiny {$t_2^{(1)}$}}
\uput{0}[0](.3,0){\Tiny {$t_1^{(1)}$}}
\uput{0}[0](2.35,2.75){ {$\zeta$}}
\uput{0}[0](-1,2){$K=\hat{\zeta}$}
\end{pspicture}
\caption{Special systems of generators}
\label{FigPlatGen}
\end{center}
\end{figure}
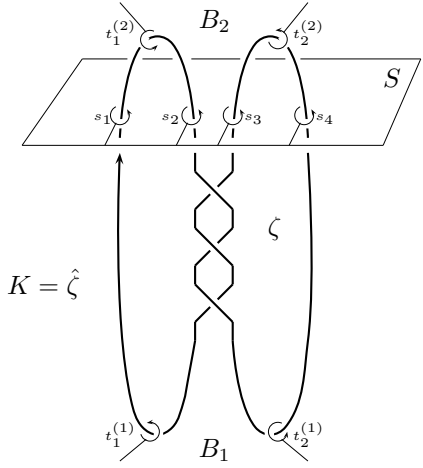

	Each knot $K \subset S^3$ can be presented as a $2n$-plat $\hat{\zeta}$. Here $\hat{\zeta}$ is obtained from the $2n$-braid $\zeta \in B_{2n}$ by closing it with $2n$ half circles as on  Fig.~\ref{FigPlatGen}. 
	Such a presentation of $K$ as plat gives rise to a {splitting} of its exterior of the form $M_K = B_1 \cup_S B_2$, where $B_1, B_2$ are two handlebodies of genus $n$ and $S = B_1 \cap B_2 = S^2 \setminus N(K)$ is a $2n$-punctured $2$-sphere (\cf~\cite[Section 3]{Heu:2003} and~\cite{VolumeForm}). This decomposition is similar to the Heegaard splitting used in the construction of the Casson invariant. It also gives rise to {special systems of generators} for $\pi_1 (B_i)$, $i=1,2$,  and for $\pi_1 (S)$ respectively denoted $\mathcal{T}_i = \{t_j^{(i)}, 1\leqslant j \leqslant n\}$ and $\mathcal{S} = \{s_j, 1 \leqslant j\leqslant 2n\}$ (see Fig.~\ref{FigPlatGen}). In fact, these systems depend on the orientation of $S^3$.

	The group $\pi_1(B_i)$ is the free group with basis $\mathcal{T}_i = \{t_j^{(i)} \; |\; 1 \leqslant j \leqslant n\}$, $i=1, 2$. The group $\pi_1(S)$ admits the finitely presentation $\pi_1(S) = \langle s_1, \ldots, s_{2n} \; | \; s_1 \cdots s_{2n}\rangle$. Furthermore, each element of $\mathcal{T}_i$ and $\mathcal{S}$ is a {meridian} of $K$. In particular, all these elements are {conjugate} in $G_K$. This obvious remark will be essential in the sequel.

	The inclusions $S \hookrightarrow B_i$ and $B_i \hookrightarrow M_K$, $i=1,2$, give rise to the commutative diagram 
\begin{equation}
\xymatrix@ur{
\pi_1 (S) \ar[r]^-{\kappa_1} \ar[d]_-{\kappa_2} 
							& \pi_1 (B_1) \ar[d]^-{p_1} \\
\pi_1 (B_2) \ar[r]_-{p_2} & \pi_1 (M_K)}\!\!=G_K
\label{E:diagramGroups}
\end{equation}
Each homomorphism of the previous diagram is onto.
The Seifert-Van Kampen Theorem and diagram~(\ref{E:diagramGroups}) combine to yield the following presentation of $G_K$: 
\begin{equation}\label{EQ:Presentation}
G_K=  \langle t^{(1)}_i, t^{(2)}_i ,\; 1\leqslant i \leqslant n \;|\;  p_1 \circ \kappa_1 (s_j) =  p_2 \circ \kappa_2 (s_j), \; 1\leqslant j \leqslant 2n-1 \rangle.
\end{equation}
Observe that~(\ref{EQ:Presentation}) is a particular {Wirtinger presentation} of $G_K$.

\subsection{Representation subspaces}
	Corresponding to a plat decomposition of $K$---which gives rise to the splitting $M_K = B_1 \cup_S B_2$---we introduce certain special re\-pre\-sen\-ta\-tion subspaces of $R(B_i)$, $i=1,2$, and $R(S)$, see~\cite[Section 3]{Heu:2003} and~\cite{VolumeForm}.
	
	Consider a representation $\rho$ of $G_K$ in $\SU$ and look at its restrictions $\rho_i = \rho \circ p_i$ and $\rho_S = \rho \circ p_i \circ \kappa_i$. In this way we do not obtain all the representations of $\pi_1(B_i)$ or $\pi_1(S)$ because {all} the generators of these groups are conjugate to the meridian of $K$. For this reason we not only consider groups but marked groups. A marked group is a pair $(G, \mathcal{G})$ where $G$ is a finitely generated group and $\mathcal{G}$ is a fixed finite set of generators of $G$. 
	For the marked group $(G, \mathcal{G})$, let $R^{\mathcal{G}}(G)$ denote the subset of $R(G) \setminus C(G)$ defined by
\begin{equation}\label{DefR}
R^{\mathcal{G}}(G)=\{\rho \in R(G)\setminus C(G) \;|\; \mathrm{Tr}{\rho(s)}=\mathrm{Tr}{\rho(t)} \, \forall s, t \in \mathcal{G}\}.
\end{equation}
It is obvious that the action by conjugation of $\SO$ leaves $R^{\mathcal{G}}(G)$ invariant. Thus we set $\widehat{R}^{\mathcal{G}}(G)=\widetilde{R}^{\mathcal{G}}(G)/\SO$ where $\widetilde{R}^{\mathcal{G}}(G)=R^{\mathcal{G}}(G)\cap \widetilde{R}(G)$. We observe that $\widetilde{R}^{\mathcal{G}}(G)$ can be identified with the total space of a principal $\SO$-bundle with  $\widehat{R}^{\mathcal{G}}(G)$ as base space. 

\begin{notation*}
We write $\widehat{R}^{\mathcal{T}_i}(B_i) = \widehat{R}^{\mathcal{T}_i}(\pi_1(B_i))$, $i = 1, 2$, $\widehat{R}^{\mathcal{S}}(S) = \widehat{R}^{\mathcal{S}}(\pi_1(S))$ etc. 
\end{notation*}	

	M. Heusener proved in~\cite{Heu:2003} that  $\widehat{R}^{\mathcal{T}_i}(B_i)$ is a $(2n-2)$-manifold and that $\widehat{R}^{\mathcal{S}}(S)$ is a $(4n-5)$-manifold. 	
	Corresponding to diagram~(\ref{E:diagramGroups}) we obtain the commutative diagram
\begin{equation}
\xymatrix@ur{
\widehat{R}^\mathcal{S}(S) 
	&\widehat{R}^{\mathcal{T}_1}(B_1) \ar[l]_-{\widehat{\kappa}_1}\\
\widehat{R}^{\mathcal{T}_2}(B_2) \ar[u]^-{\widehat{\kappa}_2}
		& \widehat{R}(M_K)  \ar[u]_-{\widehat{p}_1} \ar[l]^-{\widehat{p}_2}}
\label{E:diagramRepsQ}
\end{equation}
In diagram~(\ref{E:diagramRepsQ}) all arrows are inclusions. Therefore we can see $\widehat{R}(M_K)$ as the intersection of the images of $\widehat{R}^{\mathcal{T}_1}(B_1)$ and $\widehat{R}^{\mathcal{T}_2}(B_2)$ inside $\widehat{R}^{\mathcal{S}}(S)$. The commutative diagram~(\ref{E:diagramRepsQ}) is the main ingredient to define (generically) a volume form on the $\SU$-representation space of $G_K$.

	We can prove that a representation $\rho$ of $G_K$ in $\SU$ is {regular} if and only if the images in $\widehat{R}^{\mathcal{S}}(S)$ of the manifolds $\widehat{R}^{\mathcal{T}_1}(B_1)$ and $\widehat{R}^{\mathcal{T}_2}(B_2)$ intersect {transversally} at $[\rho]$ (see~\cite[Proposition 3.3]{Heu:2003}).
	
\subsection{Construction of the volume form}\label{SS:Construction}
	The construction of our volume form on $\Reg{K}$ is based on the two following facts:
\begin{itemize}
  \item If $\rho$ is a regular representation, then 
\begin{equation}\label{E:SuiteExacte}
\xymatrix@1@-.5pc{0 \ar[r] & \tangent{[\rho]}{\widehat{R}(M_K)} \ar[r]^-{D_{[\rho]}\widehat{p}} & \tangent{[\rho_1]}{\widehat{R}^{\mathcal{T}_1}(B_1)} \oplus \tangent{[\rho_2]}{\widehat{R}^{\mathcal{T}_2}(B_2)} \ar[r]^-{D_{[\rho]}\widehat{\kappa}} & \tangent{[\rho_S]}{\widehat{R}^{\mathcal{S}}(S)} \ar[r] & 0,}
\end{equation}
is a {short exact sequence} (cf. diagram~(\ref{E:diagramRepsQ}) and the fact that $\widehat{R}^{\mathcal{T}_1}(B_1)$ and $\widehat{R}^{\mathcal{T}_2}(B_2)$ intersect transversally et $[\rho]$).
  \item The representation subspace $\widehat{R}^{\mathcal{T}_i}(B_i)$ (resp. $\widehat{R}^{\mathcal{S}}(S)$) has a natural {$(2n-2)$-volume form} $v^{{\widehat{R}^{\mathcal{T}_i}(B_i)}}$ (resp. a natural {$(4n-5)$-volume form} $v^{\widehat{R}^{\mathcal{S}}(S)}$). The constructions of these volume forms are postponed to Subsection~\ref{SS:FormeVolNat}.
\end{itemize}
If we assume these two facts, then the $1$-volume form $\vol{\Hat{\zeta}}_{[\rho]}$ is defined using the notation of Subsection~\ref{S:Vol} and the exactness of~(\ref{E:SuiteExacte}) by 
\begin{equation}\label{E:DefVol}
\vol{\Hat{\zeta}}_{[\rho]} = (-1)^n \left( v^{{\widehat{R}^{\mathcal{T}_1}(B_1)}}_{[\rho_1]} \wedge v^{{\widehat{R}^{\mathcal{T}_2}(B_2)}}_{[\rho_2]}\right) / v^{\widehat{R}^{\mathcal{S}}(S)}_{[\rho_S]}.
\end{equation}
	In this way we locally construct a $1$-volume form $\vol{\hat{\zeta}} : [\rho] \mapsto \vol{\Hat{\zeta}}_{[\rho]} $ on the $1$-dimensional manifold $\Reg{\hat{\zeta}}$. 
	
	As it stands, $\omega^{\hat{\zeta}}$ is defined in terms of the plat decomposition of the knot and apparently depends on it. 
	In fact, we can prove that $\vol{\hat{\zeta}}$ is independent on $\hat{\zeta}$ and is moreover an invariant of knots, but this result will not be essential in the sequel. The proof of the invariance is based on a theorem of Birman-Reidemeister, which is similar for the plats to the Markov Theorem for the closed braids (see~\cite{VolumeForm} or the author  Ph.D.~\cite{JDTHESE} for details). Further observe that it is precisely the sign $(-1)^n$ which ensures the invariance of $\omega^{\hat{\zeta}}$. Thus $\Reg{K}$ is endowed with a well-defined $1$-volume form denoted $\omega^K$ which is constructed ``à la Casson". 

\begin{remark}\label{RemOr}
	The orientation induced by $\vol{K}$ is precisely the orientation on $\Reg{K}$ defined in~\cite{Heu:2003}, see~\cite{VolumeForm}.
\end{remark}

\subsection{Volume forms on the representation subspaces}
\label{SS:FormeVolNat}

The construction is based on the following steps (see~\cite{VolumeForm} for more details).
\begin{itemize}
  \item The Lie group $\SU$ is endowed with the $3$-volume form $\eta$ induced by the basis $\{\ii, \jj, \kk\}$. Similarly $\eta$ denotes the $3$-volume form on $\SO = \SU/\{\pm \I\}$ deduced from the one of $\SU$. Considering the fact that the trace-function $\mathrm{Tr} : \SU\setminus\{\pm \I\} \to (-2,2)$ is a submersion, we define a $2$-volume form $\nu$ on the $2$-sphere $S^2=\{A \in \SU \;|\; \mathrm{Tr}(A)=0\}$ exploiting the ``base $\wedge$ fiber" condition. Explicitly, for all $A \in S^2$, we have the short exact sequence
\[
\xymatrix@1@-.6pc{0 \ar[r] &\tangent{A}{S^2} \ar[r] & \tangent{A}{\SU} \ar[r] & \tangent{0}{(-2,2)} \ar[r] & 0,}
\] 
in which $\SU$ is endowed with the $3$-volume form induced by $\{\ii, \jj, \kk\}$ and $(-2, 2)$ is endowed with the usual $1$-volume form. The $2$-volume form $\nu$ on $S^2$ is the unique compatible volume form with the two others. 
  \item The map $R(B_i) \to \SU^n$ defined by $\rho \mapsto \left( {\rho(t^{(i)}_1), \ldots, \rho(t^{(i)}_n)} \right)$ is an isomorphism which allows us to identify $R(B_i)$ with $\SU^n$. Using the inclusion
$$(-2, 2) \times (S^2)^n   \to   \SU^n,\; (2\cos(\theta), P_1, \ldots, P_n)  \mapsto  (\cos(\theta) + \sin(\theta)P_i)_{1\leqslant i\leqslant n},$$ we identify $R^{\mathcal{T}_i}(B_i)$ with the product $(-2, 2) \times (S^2)^n$. As a consequence, the $(2n+1)$-manifold $R^{\mathcal{T}_i}(B_i)$ is endowed with a natural $(2n+1)$-volume form $v^{R^{\mathcal{T}_i}(B_i)}$, namely the one induced by the product volume form on $(-2, 2) \times (S^2)^n$.

	Next, we define a $(4n-2)$-volume form on $R^{\mathcal{S}}(S)$ as follows. Let $D^*$ be the $2n$-punctured disk $S \setminus \{\infty\}$. The fundamental group of $D^*$ is the free group of rank $2n$ generated by $\mathcal{S} = \{s_1, \ldots, s_{2n}\}$ (see Fig.~\ref{FigPlatGen}). Let  $U$ be the subgroup normally generated by the product $s_1 \cdots s_{2n}$, we have $\pi_1 (S) = \pi_1 (D^*) /U$. In~\cite[Lemma 3.1]{Heu:2003}, Heusener proved that the map $\psi : R^{\mathcal{S}}(D^*) \to \SU$ defined by $\psi(\rho) = \rho(s_1\cdots s_{2n})$ is surjective, that $R^{\mathcal{S}}(S)=\psi^{-1}(\I)$ and that the set of critical points of $\psi$ coincides with  the set of abelian $\SU$-representations of $\pi_1 (D^*)$.
	Thus, we have the short exact sequence
$$
\xymatrix@1@-.6pc{0 \ar[r] & \tangent{\rho}{\widetilde{R}^{\mathcal{S}}(S)} \ar[r] & \tangent{\rho}{\widetilde{R}^{\mathcal{S}}(D^*)} \ar[r] & \su \ar[r] & 0.}
$$
In this sequence, $\su$ is endowed with the $3$-volume form induced by $\{\ii, \jj, \kk\}$ and $\widetilde{R}^{\mathcal{S}}(D^*)$ is endowed with the natural $(4n+1)$-volume form (because $\widetilde{R}^{\mathcal{S}}(D^*)$  is an open subset of $R^{\mathcal{S}}(D^*) \cong (-2,2) \times (S^2)^{2n}$ which is endowed with the product volume form).  
	Then $\widetilde{R}^{\mathcal{S}}(S)$ is endowed with the unique $(4n-2)$-volume form $v^{\widetilde{R}^{\mathcal{S}}(S)}$ which is compatible with the two others.
  \item Finally, if $(G, \mathcal{G})$ is one of the marked groups $(\pi_1 (B_1), \mathcal{T}_1), (\pi_1 (B_2), \mathcal{T}_2)$ or $(\pi_1 (S), \mathcal{S})$, then $\widetilde{R}^{\mathcal{G}}(G) \to \widehat{R}^{\mathcal{G}}(G)$ is a submersion. So we define a volume form on $\widehat{R}^{\mathcal{G}}(G)$ exploiting the ``base $\wedge$ fiber" condition. 
\end{itemize}

\section{Comparison}
\label{Comparison}
	
	This section is devoted to the comparison of the two constructions we have done in Sections~\ref{Torsion} and~\ref{Vol}. We explain how we can interpret the volume form $\vol{K}$ on $\Reg{K}$ in terms of the  Reidemeister torsion form $\tau^K$ and finally prove that they produce the same topological invariant.
	
\subsection{Comparison of torsion and volume forms}
	The main result of the article is
\begin{theorem}\label{T:Volume=Torsion}
	Let $K$ be a knot in $S^3$. For any regular representation $\rho: G_K \to \SU$ we have
		$$ {\tau}^K_{\rho} = \vol{K}_{[\rho]}.$$
\end{theorem}

	The proof of Theorem~\ref{T:Volume=Torsion} needs several steps (and is postponed to Sections~\ref{S:Volume=Torsion1} and~\ref{S:Volume=Torsion2}). The key idea to compare the two volume forms $\tau^K_\rho$ and $\vol{K}_{[\rho]}$ on $\tangent{[\rho]}{\widehat{R}(M_K)}$ is to interpret them as one and the same  Reidemeister torsion of an auxiliary chain complex. We shall discuss this point in Subsection~\ref{IdeasProof}.

\subsection{Some remarks and consequences}
\begin{enumerate}
  \item According to Remark~\ref{RemOr} the induced orientation on $\Reg{K}$ deduced from the volume form $\tau^K$ is exactly the one defined on $\Reg{K}$ by Heusener in~\cite{Heu:2003}.
  \item The Casson-Lin invariant $\lambda_{CL}(K)$, defined in~\cite{Lin:1992}, 
 can be interpreted via the sign-determined Reidemeister torsion as follows (\cf~\cite[Proposition 5.4 \& Lemma 5.5]{Heu:2003}).  Let $(G, \mathcal{G})$ be a marked group. 
 Consider $R^{\mathcal{G}}_{0}(G) = \{\rho \in R^{\mathcal{G}} \, |\, \mathrm{Tr}(\rho(t)) = 0 \; \forall t \in \mathcal{G}\}$ and $\widehat{R}^{\mathcal{G}}_{0}(G) = \widetilde{R}^{\mathcal{G}}_{0}(G)/\SO$ (with $\widetilde{R}^{\mathcal{G}}_{0}(G) = R^{\mathcal{G}}_{0}(G) \cap \widetilde{R}(G)$). Set $\widehat{Q}_{i}^{(0)} = \widehat{\kappa}_{i}(\widehat{R}^{\mathcal{T}_{i}}_{0}(B_{i}))$ for $i=1, 2$. 
  The intersection $\widehat{Q}_{1}^{(0)} \cap \widehat{Q}_{2}^{(0)}$ is compact  and the Casson-Lin invariant of the closed braid $\sigma^{\wedge}$, $\sigma \in B_{n}$, is the intersection number: $$\lambda_{CL}(\sigma^{\wedge}) = (-1)^{n}\left\langle {\widehat{Q}_{1}^{(0)}, \widehat{Q}_{2}^{(0)}}\right\rangle_{\widehat{R}^{\mathcal{S}}_{0}(S)}$$ see~\cite[Theorem 18]{Lin:1992} and~\cite[Section 5.3]{Heu:2003}. Keep in mind the notation of Subsection~\ref{RegularRep}. If {all} the representations in $\widehat{Q}_{1}^{(0)} \cap \widehat{Q}_{2}^{(0)}$ are $m$-regular, \ie if $\widehat{Q}_{1}^{(0)}$ and $\widehat{Q}_{2}^{(0)}$ intersect {transversally}, then
\[
\lambda_{CL}(\sigma^{\wedge}) = \sum_{[\rho] \in \widehat{Q}_{1}^{(0)} \cap \widehat{Q}_{2}^{(0)}} \mathrm{sgn}\left( \tau^{\sigma^{\wedge}}_\rho(h^{(1)}_{\rho}(m))\right).
\]
\end{enumerate}

\subsection{Main ideas of the proof}
\label{IdeasProof}

	Let $K$ be an oriented knot in $S^3$ presented as a $2n$-plat $\hat{\zeta}$ as in Section~\ref{Vol} and let $\rho : G_K \to \SU$ be a regular representation. The plat presentation $\hat{\zeta}$ of $K$ induces the Wirtinger presentation~(\ref{EQ:Presentation}) of $G_K$.
	
	We want to compare the Reidemeister torsion of the acyclic chain complex~(\ref{E:SuiteExacte})
\[
\mathscr{T}^\rho_* = \xymatrix@1@-.5pc{0 \ar[r] & \tangent{[\rho]}{\widehat{R}(M_K)} \ar[r] & \tangent{[\rho_1]}{\widehat{R}^{\mathcal{T}_1}(B_1)} \oplus \tangent{[\rho_2]}{\widehat{R}^{\mathcal{T}_2}(B_2)} \ar[r] & \tangent{[\rho_S]}{\widehat{R}^{\mathcal{S}}(S)} \ar[r] & 0,}
\] 
which allows us to define the volume form $\omega^K_{[\rho]}$, to the Reidemeister torsion of the non-acyclic $(Ad \circ \rho)$-twisted chain complex~(\ref{twistX}) of the knot exterior
$$\mathscr{X}^\rho_* = \xymatrix@1@-.5pc{
0 \ar[r] & \su \ar[r]^-{d^{\rho}_2} & \su^{2n} \ar[r]^-{d^{\rho}_1} & \su^{2n-1} \ar[r] & 0.
}$$
For this purpose we introduce and use the \emph{auxiliary chain complex}
\[
\mathscr{D}^\rho_* =  \xymatrix@1@-.7pc{0 \ar[r] & \Der{\rho}{(M_K)} \ar[r]^-{p^*} & \Der{\rho_1}{(B_1)} \oplus \Der{\rho_2}{(B_2)} \ar[r]^-{\kappa^*} & \Der{\rho_S}{(S)} \ar[r] & 0.}
\]
Applying Propositions~\ref{isoRandH} and~\ref{IsomTH} we get $\Der{\rho}{(G_K)} \cong \tangent{\rho}{\widetilde{R}(M_K)}$, $\Der{\rho_i}{(\pi_1(B_i))} \cong \tangent{\rho_i}{\widetilde{R}(B_i)}$ and $\Der{\rho_S}{(\pi_1(S))} \cong \tangent{\rho_S}{\widetilde{R}(S)}$. 
The chain complex $\mathscr{D}^\rho_*$ is {not acyclic} and can be consider as a sort of ``Mayer-Vietoris sequence". The advantage in using $\mathscr{D}^\rho_*$ belongs to the fact that, in opposition to $\mathscr{X}^\rho_*$, its homology in dimension $1$ vanishes and it has just homology in dimension $0$: $H_0(\mathscr{D}^\rho_*) \cong H^2_\rho(M_K) \cong H_0(\mathscr{X}^\rho_*)$.  
We will discuss  the chain complex $\mathscr{D}^\rho_*$ in detail in Section~\ref{AuxComplex}.

The main ingredient to perform the comparison of the torsions is the {multipli\-ca\-ti\-vi\-ty} of the sign-determined Reidemeister torsion. Let us recall this property.
	
\begin{MultLemma}\label{SS:Multiplicativite}	Let 
\begin{equation}\label{EQM}
\xymatrix@1@-.6pc{0 \ar[r] & C'_* \ar[r] & C^{}_* \ar[r] & C_*'' \ar[r] & 0}
\end{equation}
be an exact sequence of chain complexes. For each $i$, fix volume forms on  $C_i'$, $C^{}_i$, $C_i''$, $H_i(C'_*)$, $H_i(C^{}_*)$ and $H_i(C''_*)$. Associated to~(\ref{EQM}) is the long exact sequence in homology
\[
\xymatrix@1@-.6pc{ \cdots \ar[r] & H_i(C'_*) \ar[r] & H_i(C^{}_*) \ar[r] & H_i(C''_*) \ar[r] & H_{i-1}(C'_*) \ar[r] & \cdots.}
\] 
Let $\mathcal{H}_*$ denote this acyclic chain complex and endow $\mathcal{H}_{3i+2} = H_i(C'_*)$, $\mathcal{H}_{3i+1} = H_i(C^{}_*)$ and $\mathcal{H}_{3i} = H_i(C''_*)$ with the volume form on  $H_i(C'_*)$, $H_i(C^{}_*)$ and $H_i(C''_*)$ respectively.

If, for all $i$, the volume forms on the components of the short exact sequence $\xymatrix@1@-.6pc{0 \ar[r] & C'_i \ar[r] & C^{}_i \ar[r] & C_i'' \ar[r]& 0}$ are compatible, then
	$$\mathrm{Tor}(C^{}_*) = (-1)^{\alpha(C'_*, C''_*) + \varepsilon(C'_*, C^{}_*, C''_*)} \; \mathrm{Tor}(C'_*) \cdot \mathrm{Tor}(C''_*) \cdot \mathrm{tor}(\mathcal{H}_*),$$
where $$\alpha(C'_*, C''_*) = \sum_{i \geqslant 0} \alpha_{i-1}(C'_*) \alpha_i(C''_*) \in \ZZ/2\ZZ$$
and $$\varepsilon(C'_*, C^{}_*, C''_*) = \sum_{i \geqslant 0} [(\beta_i(C^{}_*)+1)(\beta_i(C_*')+\beta_i(C_*'')) + \beta_{i-1}(C_*')\beta_i(C_*'')] \in \ZZ/2\ZZ.$$
\end{MultLemma}
	The proof of the Multiplicativity Lemma is a careful computation based on linear algebra, see~\cite[Lemma 3.4.2]{Turaev:1986} and~\cite[Theorem 3.2]{Milnor:1966}.

	Keeping in mind this property, the proof of Theorem~\ref{T:Volume=Torsion} is done in Sections~\ref{S:Volume=Torsion1} and~\ref{S:Volume=Torsion2}\,; it is quite natural and essentially consists of two parts. Firstly (see Section~\ref{S:Volume=Torsion1}) we compare the sign-determined Reidemeister torsions of the chain complexes $\mathscr{T}^\rho_*$ and $\mathscr{D}^\rho_*$. Secondly (see Section~\ref{S:Volume=Torsion2}) we compare the sign-determined Reidemeister torsions of the chain complexes $\mathscr{D}^\rho_*$ and $\mathscr{X}^\rho_*$ to obtain the required equality between the torsion form $\tau^K$ and the volume form $\omega^K$.

\section{The auxiliary chain complex}
\label{AuxComplex}
 
	The proof of Theorem~\ref{T:Volume=Torsion} uses some more properties about the notion of twisted derivations; this material is described in the following subsections. The notion of twisted derivations gives a method to identify the tangent bundles of the $\SU$-re\-pre\-sen\-ta\-tion subspaces $\widetilde{R}^{\mathcal{T}_i}(B_i)$ and $\widetilde{R}^{\mathcal{S}}(S)$. These identifications are used to effectively compare the volume forms $\vol{K}$ and $\tau^K$. 

\subsection{Volume forms on the twisted derivation spaces}\label{SS:Der}
	Suppose that $K$ is presented as a $2n$-plat and keep the notation of Section~\ref{Vol}. The plat presentation of $K$ gives the splitting $M_K = B_1 \cup_S B_2$, where $B_1, B_2$ are handlebodies and $S$ is a $2n$-punctured $2$-sphere.
	
\begin{notation*}
We write $\Der{\rho}{(M_K)} = \Der{\rho}{(G_K)}$, $\mathrm{Der}_{\rho_i}(B_i) = \mathrm{Der}_{\rho_i}(\pi_1(B_i))$ and $\mathrm{Der}_{\rho_S}(S) = \mathrm{Der}_{\rho_S}(\pi_1(S))$, here $\rho \in R(M_K)$.
\end{notation*}
	
	The choice of the system of generators $\mathcal{T}_i$ (resp. $\mathcal{S}$) for the group  $\pi_1 (B_i)$, $i=1, 2$, (resp. $\pi_1 (S)$) gives
\begin{fact}
	For $i=1, 2$, $\Der{\rho_i}{(B_i)} \cong \su^n$ and  $\Der{\rho_S}{(S)} \cong \su^{2n-1}$.
\end{fact}
\begin{proof} 
We use the fact that $\pi_1(B_i)$ (resp. $\pi_1(S)$) is the free group of rank $n$ (resp. is isomorphic to the free group of rank $2n-1$).
Explicitly the isomorphisms are
\begin{itemize}
  \item $\varphi_i : \su^n \to \mathrm{Der}_{\rho_i}(B_i)$ defined by $\varphi_i : \mathbf{x} = (x_1, \ldots, x_n) \mapsto d^{\mathbf{x}}$, where $d^{\mathbf{x}} \in \mathrm{Der}_{\rho_i}(B_i)$ is such that $d^{\mathbf{x}}(t^{(i)}_j) = x_j$, for all $j=1, \ldots, n$.
  \item  $\varphi_S : \su^{2n-1} \to \mathrm{Der}_{\rho_S}(S)$ defined by $\varphi_S : \mathbf{x} = (x_1, \ldots, x_{2n-1}) \mapsto d^{\mathbf{x}}$, where $d^{\mathbf{x}} \in \mathrm{Der}_{\rho_S}(S)$ is such that $d^{\mathbf{x}}(s_j) = x_j$, for all $j=1, \ldots, 2n-1,$ and $$d^{\mathbf{x}}(s_{2n}) = - Ad_{\rho_S(s_1 \cdots s_{2n-1})^{-1}} (x_1) - \cdots - Ad_{\rho_S(s_{2n-1})^{-1}} (x_{2n-1}).$$
\end{itemize}
\end{proof}

	The isomorphism $\mathrm{Der}_{\rho_i}(B_i) \cong \su^n$ allows us to define a {volume form} on $\mathrm{Der}_{\rho_i}(B_i)$ namely the one induced by the product volume form on $\su^n$. Similarly, the isomorphism $\Der{\rho_S}{(S)} \cong \su^{2n-1}$ allows us to define a {volume form} on $\mathrm{Der}_{\rho_S}(S)$ in the same way. But it is possible to define an \emph{a priori} different volume form on $\mathrm{Der}_{\rho_S}(S)$ in another ``natural" way.
	Consider the $2n$-punctured disc $D^* = S \setminus \{\infty\}$ and the map $\psi : \mathrm{Der}_{\rho_S}(D^*) \to \su$ defined by $\psi(d) = d(s_1 \cdots s_{2n})$. Thus, we have the short exact sequence
\begin{equation}\label{S}
\xymatrix@1@-.7pc{0 \ar[r] & \mathrm{Der}_{\rho_S}(S) \ar[r] & \mathrm{Der}_{\rho_S}(D^*) \ar[r]^-{\psi} & \su \ar[r] & 0.}
\end{equation}
Recall that $\mathrm{Der}_{\rho_S}(D^*) \cong \su^{2n}$ is endowed with the volume form induced by the product volume form.
\begin{claim}
The volume forms deduced from the product volume forms on $\su^m$ on each component of the exact sequence~(\ref{S}) are compatible.
\end{claim}
\begin{proof}
Let $d^x \in \mathrm{Der}_{\rho_S}(D^*)$ be defined by $d^x(s_j) = 0$, for all $j=1, \ldots, 2n-1$, and $d^x(s_{2n}) = Ad_{\rho_S(s_{2n})}(x)$. 
We consider the section $s : \su \to \mathrm{Der}_{\rho_S}(D^*)$ of $\psi$ defined by $s(x) = d^x$. 
	The previous identifications for $\mathrm{Der}_{\rho_S}(S)$ and $\mathrm{Der}_{\rho_S}(D^*)$ gives us the commutative diagram
\[
\xymatrix@-.7pc{0 \ar[r] & \mathrm{Der}_{\rho_S}(S) \ar[r] \ar[d]^-{\cong} & \mathrm{Der}_{\rho_S}(D^*) \ar[r]^-{\psi} \ar[d]^-{\cong} & \su \ar[r] \ar[d]^-{=} & 0 \\
0 \ar[r] & \su^{2n-1} \ar[r]^-i & \su^{2n} \ar[r]^-j & \su \ar[r] \ar@/^1pc/[l]^s & 0}
\]
in which
\begin{multline}
i(x_1, \ldots, x_{2n-1}) = \\ \left( x_1, \ldots, x_{2n-1}, - Ad_{\rho_S(s_1 \cdots s_{2n-1})^{-1}} (x_1) - \cdots - Ad_{\rho_S(s_{2n-1})^{-1}} (x_{2n-1}) \right) \notag
\end{multline}
and
$$j(x_1, \ldots, x_{2n}) = x_1 + Ad_{\rho_S(s_1)}(d(s_2)) + \cdots + Ad_{\rho_S(s_1 \cdots s_{2n-1})}(d(s_{2n})).$$
Moreover the section $s$ of $\psi$ can be identified with the one of $j$ also denoted $s : \su \to \su^{2n}$ and defined by $s(x) = (0, \ldots, 0, Ad_{\rho_S(s_{2n})}(x))$.

	The determinant of the transition matrix from the canonical product basis $$\{\ii_1, \jj_1, \kk_1, \ldots, \ii_{2n}, \jj_{2n}, \kk_{2n}\}$$ of $\su^{2n}$ to the basis $$\left\{ {i(\ii_1), i(\jj_1), i(\kk_1), \ldots, i(\ii_{2n-1}), i(\jj_{2n-1}), i(\kk_{2n-1}), s(\ii), s(\jj), s(\kk)} \right\}$$ of $\su^{2n}$ is $+1$. So, the volume form on $\mathrm{Der}_{\rho_S}(S)$ which is compatible with the two others on the components of the exact sequence~(\ref{S}) can be identified with the product volume form on $\su^{2n-1}$.
\end{proof}

\subsection{The auxiliary derivation chain complex}
\label{AuxDerCW}
	
	The Mayer-Vietoris sequence for the twis\-ted cohomology associated to the splitting $M_K = B_1 \cup_S B_2$ and to the irreducible representation $\rho : G_K \to \SU$ reduces to
\[
\xymatrix@1@-.7pc{0 \ar[r] & H^1_\rho(M_K) \ar[r] & H^1_{\rho_1}(B_1) \oplus H^1_{\rho_2}(B_2) \ar[r] & H^1_{\rho_S}(S) \ar[r] & H^2_\rho(M_K) \ar[r] & 0.}
\]
From the splitting $M_K = B_1 \cup_S B_2$ we deduce the \emph{auxiliary derivation chain complex}
\[
\mathscr{D}^\rho_* = \xymatrix@1@-.7pc{0 \ar[r] & \Der{\rho}{(M_K)} \ar[r]^-{p^*} & \Der{\rho_1}{(B_1)} \oplus \Der{\rho_2}{(B_2)} \ar[r]^-{\kappa^*} & \Der{\rho_S}{(S)} \ar[r] & 0.}
\]
Here
\begin{align*}
	p^* &:   \Der{\rho}{(M_K)} \to \Der{\rho_1}{(B_1)} \oplus \Der{\rho_2}{(B_2)}, \ d \mapsto (d \circ p_1, -d \circ p_2) \\ 
\text{and} \;	\kappa^* &:  \Der{\rho_1}{(B_1)} \oplus \Der{\rho_2}{(B_2)} \to \Der{\rho_S}{(S)}, \; (d_1, d_2) \mapsto d_1 \circ \kappa_1 + d_2 \circ \kappa_2.
\end{align*}

The equalities $\ker p^* = 0$, $\im p^* = \ker \kappa^*$ and $\coker \kappa^* = H^2_\rho(M_K)$ give us
\begin{fact}\label{HofD}
	The homology of $\mathscr{D}^\rho_*$ reduces to $H_*(\mathscr{D}^\rho_*) = H_0(\mathscr{D}^\rho_*) \cong H^2_\rho(M_K)$.
\end{fact}

	We now describe the {reference volume form} we choose on each component of the chain complex $\mathscr{D}^\rho_*$. Fact~\ref{SS:Der} tells us that $\varphi_i : \Der{\rho_i}{(B_i)} \to \su^n$ (resp. $\varphi_S : \Der{\rho_S}{(S)} \to \su^{2n-1}$) is an isomorphism. So $\Der{\rho_i}{(B_i)}$ (resp. $\Der{\rho_S}{(S)}$) is endowed with the volume form induced by $\varphi_i$ (resp. $\varphi_S$) and by the product volume form on $\su^n$  (resp. $\su^{2n-1}$). We endow $\Der{\rho_1}{(B_1)} \oplus \Der{\rho_2}{(B_2)}$ with the product volume form corrected by the factor $(-1)^n$ (as in equation~(\ref{E:DefVol})).
	To produce a reference volume form on $\Der{\rho}{(M_K)}$ we use the exactness of sequence~(\ref{ExSH})
\[
\xymatrix@1@-.6pc{0 \ar[r] & \su \ar[r] & \Der{\rho}{(M_K)} \ar[r] & H^1_\rho(M_K) \ar[r] & 0.}
\]
Here $\su$ is endowed with the $3$-volume form $\ii \wedge \jj \wedge \kk$ and $H^1_\rho(M_K) \cong \tangent{[\rho]}{\widehat{R}(M_K)}$ is endowed with the $1$-volume form  $\omega^K_{[\rho]}$. Using compatibility we endow $\Der{\rho}{(M_K)}$ with the $4$-volume form $w_\rho$ such that $w_\rho/\omega^K_{[\rho]} = \ii \wedge \jj \wedge \kk$ (we use the notation described in Subsection~\ref{S:Vol}).

	The complex $\mathscr{D}^\rho_*$ is based using the convention of Subsection~\ref{AlgPre}; that is, we endow each $\mathscr{D}^\rho_i$ with any basis which has volume one. 

\subsection{Tangent bundles of $\widetilde{R}^{\mathcal{T}_i}(B_i)$, $i=1, 2$, and $\widetilde{R}^{\mathcal{S}}(S)$}
\label{Tspace}
	
	Let $(G, \mathcal{G})$ be a marked group. Set $\mathcal{G} = \{g_1, \ldots, g_m\}$ and fix $\rho \in \widetilde{R}^\mathcal{G}(G)$. Just as we are interested in the representation subspace $R^{\mathcal{G}}(G)$ we introduce the subspace 
$$\mathrm{Der}^{\mathcal{G}}_{\rho} (G)=\{d \in \Der{\rho}{(G)} \; | \; \langle d(g_j), P^{\rho}(g_j) \rangle = \langle d(g_k), P^{\rho}(g_k) \rangle, 1\leqslant j, k \leqslant m \}$$
of $\Der{\rho}{(G)}$. We also consider the subspace
$${\mathrm{Der}^{\mathcal{G}}_{\rho}(G)}_0=\{d \in \Der{\rho}{(G)} \; | \; \langle d(g_j), P^{\rho}(g_j) \rangle = 0, 1\leqslant j \leqslant m \}$$
 of $\mathrm{Der}^{\mathcal{G}}_{\rho}(G)$.
The equalities between the scalar products can be considered as an ``infinitesimal" counterpart of the equalities between traces in~(\ref{DefR}).

\begin{notation*}
	We write $\mathrm{Der}^{\mathcal{T}_i}_{\rho_i}(B_i) = \mathrm{Der}^{\mathcal{T}_i}_{\rho_i}(\pi_1(B_i))$ and  $\mathrm{Der}^{\mathcal{S}}_{\rho_S}(S) = \mathrm{Der}^{\mathcal{S}}_{\rho_S}(\pi_1(S))$ etc.
\end{notation*}

	The sequel of this subsection is devoted to a careful analysis of the subspaces $\mathrm{Der}^{\mathcal{T}_i}_{\rho_i}(B_i)$ and $\mathrm{Der}^{\mathcal{S}}_{\rho_S}(S)$ where $\rho \in \widetilde{R}(M_K)$. In fact, we shall prove the two isomorphisms   $\mathrm{Der}^{\mathcal{T}_i}_{\rho_i}(B_i) \cong \tangent{\rho_i}{\widetilde{R}^{\mathcal{T}_i}(B_i)}$ and $\mathrm{Der}^{\mathcal{S}}_{\rho_S}(S) \cong \tangent{\rho_S}{\widetilde{R}^{\mathcal{S}}(S)}$.
	Let us begin with a lemma.
\begin{lemma}\label{MainLemma}
We have $\Inn{\rho_i}{(B_i)} \subset {\mathrm{Der}^{\mathcal{T}_i}_{\rho_i}(B_i)}_0$, for $i=1, 2,$ and
 $\Inn{\rho_S}{(S)} \subset {\mathrm{Der}^{\mathcal{S}}_{\rho_S}(S)}_0$. Moreover if $d \in \Der{\rho}{(M_K)}$, then $d_i = d \circ p_i \in \mathrm{Der}^{\mathcal{T}_i}_{\rho_i} (B_i)$ for $i=1, 2$.
\end{lemma}
\begin{proof}
	The proof of this lemma is a straightforward consequence of the well-known identity $\langle x, P \rangle = \langle Ad_A (x), P\rangle,$ where $A \in \SU \setminus \{\pm \I\}$, $Ad_A$ is a rotation which fixes $P \in S^2$ and $x \in \su$. In fact, we prove the following more general result: if $G$ is a group, if $g \in G$ and if $\rho \in \widetilde{R}(G)$ is such that $\rho(g) \ne \pm \I$, then we have the properties
\begin{enumerate}
  \item if $\delta \in \mathrm{Inn}_\rho(G)$, then $\langle \delta(g), P^\rho(g) \rangle = 0$,
  \item if $d \in \mathrm{Der}_\rho(G)$ and $g$, $g'$ are conjugate, then $\langle d(g), P^\rho(g) \rangle = \langle d(g'), P^\rho(g') \rangle$. 
\end{enumerate}
For this purpose,
\begin{enumerate}
  \item If $\delta \in \Inn{\rho}{(G)}$, then there exists $a \in \su$ such that $\delta(g) = a - Ad_{\rho(g)}a$ for all $g \in G$; so that $\langle \delta(g), P^\rho(g) \rangle = \langle a, P^\rho(g) \rangle - \langle Ad_{\rho(g)}a, P^\rho(g) \rangle = 0$.
  \item If we write $g'$ as $g' = hgh^{-1}$, then
$$d(g') = d(h) - Ad_{\rho(g')}d(h) + Ad_{\rho(h)}d(g)\; \text{and}\; P^\rho(g') = Ad_{\rho(h)}P^\rho(g).$$
We thus have
\begin{multline}
\langle d(g'), P^\rho(g') \rangle = \langle d(h), P^\rho(g') \rangle - \langle Ad_{\rho(g')}d(h), P^\rho(g') \rangle \\ + \langle Ad_{\rho(h)}d(g), Ad_{\rho(h)}P^\rho(g) \rangle, \notag
\end{multline}
which gives $\langle d(g'), P^\rho(g') \rangle = \langle d(g), P^\rho(g) \rangle$.
\end{enumerate}

	The inclusions stated in the lemma straight follow from the first property; the second part arises from the second property (because $t^{(i)}_j$ is conjugate to $m$).
	\end{proof}

	The same arguments give us 
\begin{lemma}
If $\rho \in \widetilde{R}(M_K)$, then $\mathrm{Inn}_\rho(M_K) \subset {\mathrm{Der}_\rho(M_K)}_0 \subset \mathrm{Der}_\rho(M_K)$.
\end{lemma}
Here we write ${\mathrm{Der}_\rho(M_K)}_0 = \{d \in \mathrm{Der}_\rho(M_K)  \; | \; \langle d(m), P^\rho(m) \rangle = 0\}$.
\begin{remark}\label{remarkreg}
	Recall that $\rho$ is regular if and only if $\dim \mathrm{Der}_\rho(M_K) = 4$. Observe that for any regular representation $\rho$ we have either ${\mathrm{Der}_\rho(M_K)}_0 = \mathrm{Inn}_\rho(M_K)$ when $\rho$ is $m$-regular, or ${\mathrm{Der}_\rho(M_K)}_0 =\mathrm{Der}_\rho(M_K)$ when $\rho$ is regular but not $m$-regular.
\end{remark}

	If $(G, \mathcal{G})$ is one of the marked groups $(\pi_1 (B_1), \mathcal{T}_1)$, $(\pi_1 (B_2), \mathcal{T}_2)$, or $(\pi_1 (S), \mathcal{S})$, then inclusion~(\ref{TRintoDer}) induces an inclusion from $\tangent{\rho}{\widetilde{R}^\mathcal{G}(G)}$ into $\mathrm{Der}_{\rho}^\mathcal{G}(G)$. Using a dimensional argument we obtain 
\begin{lemma}\label{lemmaisomo}
If $\rho \in \widetilde{R}(M_K)$, then the inclusions $\tangent{\rho}{\widetilde{R}^{\mathcal{T}_i}(B_i)} \hookrightarrow \mathrm{Der}^{\mathcal{T}_i}_\rho (B_i)$, $i=1, 2$, and $\tangent{\rho}{\widetilde{R}^{\mathcal{S}}(S)} \hookrightarrow \mathrm{Der}^{\mathcal{S}}_\rho (S)$ are isomorphisms.
\end{lemma}

	The space $\mathrm{Der}^{\mathcal{T}_i}_\rho (B_i) \cong \tangent{\rho}{\widetilde{R}^{\mathcal{T}_i}(B_i)}$, $i=1, 2$, (resp. $\mathrm{Der}^{\mathcal{S}}_\rho (S) \cong \tangent{\rho}{\widetilde{R}^{\mathcal{S}}(S)}$) is endowed with the volume form induced by the natural one  $v^{\widetilde{R}^{\mathcal{T}_i}(B_i)}_{\rho_i}$ on $\tangent{\rho}{\widetilde{R}^{\mathcal{T}_i}(B_i)}$ (resp. $v^{\widetilde{R}^{\mathcal{S}}(S)}_{\rho_S}$ on $\tangent{\rho}{\widetilde{R}^{\mathcal{S}}(S)}$), see Subsection~\ref{SS:FormeVolNat}.  These isomorphisms between tangent subspaces and derivation subspaces are the main tool in the proof of Theorem~\ref{T:Volume=Torsion}. 
	
	For any regular representation $\rho$ of $G_K$ we have the acyclic chain complex
\begin{equation}\label{E:Derivations}
\mathscr{C}^\rho_* = \xymatrix@1@-.5pc{0 \ar[r] & \mathrm{Der}_\rho (M_K) \ar[r]^-{p^*} & \mathrm{Der}^{\mathcal{T}_1}_{\rho_1} (B_1) \oplus \mathrm{Der}^{\mathcal{T}_2}_{\rho_2} (B_2) \ar[r]^-{\kappa^*}  & \mathrm{Der}^{\mathcal{S}}_{\rho_S} (S) \ar[r] & 0,}
\end{equation}
and the isomorphisms $\Der{\rho}{(M_K)} \cong \tangent{\rho}{\widetilde{R}(M_K)}$, $\mathrm{Der}_{\rho_i}^{\mathcal{T}_i}(B_i) \cong \tangent{\rho_i}{\widetilde{R}^{\mathcal{T}_i}(B_i)}$, $i=1, 2$, and $\mathrm{Der}_{\rho_S}^{\mathcal{S}}(S) \cong \tangent{\rho_S}{\widetilde{R}^{\mathcal{S}}(S)}$. 
	We endow the components of $\mathscr{C}^\rho_*$ with {reference volume forms} as follows: 
\begin{itemize}
  \item $\mathscr{C}^\rho_2 = \mathrm{Der}_\rho (M_K)$ is endowed with the $4$-volume form $w_\rho$,
  \item $\mathscr{C}^\rho_1 = \mathrm{Der}^{\mathcal{T}_1}_{\rho_1} (B_1) \oplus \mathrm{Der}^{\mathcal{T}_2}_{\rho_2} (B_2)$ is endowed with the sign-corrected product volume form $(-1)^n (v^{\widetilde{R}^{\mathcal{T}_1}(B_1)}_{\rho_1} \wedge v^{\widetilde{R}^{\mathcal{T}_2}(B_2)}_{\rho_2})$ (as in equation~(\ref{E:DefVol})),
  \item $\mathscr{C}^\rho_0 = \mathrm{Der}^{\mathcal{S}}_{\rho_S} (S)$ is endowed with $v^{\widetilde{R}^{\mathcal{S}}(S)}_{\rho_S}$.
\end{itemize}	
	The complex $\mathscr{C}^\rho_*$ is based using the convention of Subsection~\ref{AlgPre}; that is, we endow each $\mathscr{C}^\rho_i$ with any basis which has volume one.

	We finally introduce the linear forms 
\begin{align}
	\pi_i &: \mathrm{Der}_{\rho_i}^{\mathcal{T}_i}(B_i) \to \IR, \; d \mapsto \langle d(t^{(i)}_1), P^{\rho_i}(t^{(i)}_1)\rangle, \; i=1, 2, \label{Pii}\\
	\pi_S &: \mathrm{Der}_{\rho_S}^{\mathcal{S}}(S) \to \IR, \; d \mapsto \langle d(s_1), P^{\rho_S}(s_1)\rangle.\label{PiS}
\end{align}
It is quite obvious that $\pi_i$ and $\pi_S$ are surjective, that ${\mathrm{Der}^{\mathcal{T}_i}_{\rho_i}(B_i)}_0 = \ker \pi_i,$ for $i=1, 2,$ and that ${\mathrm{Der}^{\mathcal{S}}_{\rho_S}(S)}_0 = \ker \pi_S$. These equalities allow us to define a natural volume form on each kernel---namely the one induced by the volume form on $\mathrm{Der}^{\mathcal{T}_i}_{\rho_i}(B_i)$ and $\mathrm{Der}^{\mathcal{S}}_{\rho_S}(S)$ respectively.

\section{Proof of Theorem~\ref{T:Volume=Torsion} : Part 1}
\label{S:Volume=Torsion1}

	We turn now to the proof of Theorem~\ref{T:Volume=Torsion}, which will be discussed in this section and the next. To prove the equality between the volume form $\omega^K$ and torsion form $\tau^K$ we interpret them as the Reidemeister torsion of the auxiliary chain complex 
$$\mathscr{D}^\rho_* = \xymatrix@1@-.7pc{0 \ar[r] & \Der{\rho}{(M_K)} \ar[r]^-{p^*} & \Der{\rho_1}{(B_1)} \oplus \Der{\rho_2}{(B_2)} \ar[r]^-{\kappa^*} & \Der{\rho_S}{(S)} \ar[r] & 0}.$$ 
	The proof splits into two parts. In the first one (this section) we explicitly compute the sign-determined Reidemeister torsion of $\mathscr{D}^\rho_*$. In the next one (Section~\ref{S:Volume=Torsion2}) we achieve the proof  by comparing the sign-determined Reidemeister torsion of $\mathscr{D}^\rho_*$ to the one of the twisted complex~(\ref{twistX}) of the knot exterior
$$\mathscr{X}^\rho_* = \xymatrix@1@-.5pc{
0 \ar[r] & \su \ar[r]^-{d^{\rho}_2} & \su^{2n} \ar[r]^-{d^{\rho}_1} & \su^{2n-1} \ar[r] & 0}.$$

\medskip

	The main purpose of the first part of the proof is to establish
\begin{lemma}\label{P:TorD}
	If $H_0(\mathscr{D}^\rho_*) \cong H^2_\rho(M_K)$ is based with the generator $h^{(2)}_\rho$ (see equation~(\ref{EQ:Defh2})), then $$\mathrm{Tor}(\mathscr{D}^\rho_*) = \mathrm{sgn}(\mathrm{Tor}(\mathfrak{X}_*)).$$
\end{lemma}
Here $\mathfrak{X}_*$ is the chain complex~(\ref{NoeudReel}).

The proof of Lemma~\ref{P:TorD} requires three steps. It is based on Lemma~\ref{lemmaisomo}.

\subsection{First Step}

The aim of the first step is to compare the sign-determined Rei\-de\-meis\-ter torsion of the acyclic based chain complex (see Subsection~\ref{Tspace}) 
$$
\mathscr{C}^{\rho}_* = \xymatrix@1@-.6pc{0 \ar[r] & \mathrm{Der}_\rho (M_K) \ar[r] & \mathrm{Der}^{\mathcal{T}_1}_{\rho_1} (B_1) \oplus \mathrm{Der}^{\mathcal{T}_2}_{\rho_2} (B_2) \ar[r] & \mathrm{Der}^{\mathcal{S}}_{\rho_S} (S) \ar[r] & 0}
$$
with the one of the acyclic chain complex
$$
\mathscr{T}^\rho_* = \xymatrix@1@-.6pc{0 \ar[r] & \tangent{[\rho]}{\widehat{R}(M_K)} \ar[r] & \tangent{[\rho_1]}{\widehat{R}^{\mathcal{T}_1}(B_1)} \oplus \tangent{[\rho_2]}{\widehat{R}^{\mathcal{T}_2}(B_2)} \ar[r] & \tangent{[\rho_S]}{\widehat{R}^{\mathcal{S}}(S)} \ar[r] & 0.}
$$
The components of this last complex are endowed with the compatible volume forms described in Subsections~\ref{SS:Construction}--\ref{SS:FormeVolNat}: 
\begin{itemize}
  \item $\mathscr{T}^\rho_0 = \tangent{[\rho_S]}{\widehat{R}^{\mathcal{S}}(S)}$ is endowed with $v^{\widehat{R}^{\mathcal{S}}(S)}_{[\rho_S]}$,
  \item $\mathscr{T}^\rho_1 = \tangent{[\rho_1]}{\widehat{R}^{\mathcal{T}_1}(B_1)} \oplus \tangent{[\rho_2]}{\widehat{R}^{\mathcal{T}_2}(B_2)}$ is endowed with the product volume form $$(-1)^n \left( v^{\widehat{R}^{\mathcal{T}_1}(B_1)}_{[\rho_1]} \wedge v^{\widehat{R}^{\mathcal{T}_2}(B_2)}_{[\rho_2]}\right),$$ 
  \item $\mathscr{T}^\rho_2 = \tangent{[\rho]}{\widehat{R}(M_K)}$ is endowed with the $1$-volume form $\vol{K}_{[\rho]}$ defined by equation~(\ref{E:DefVol}). 
\end{itemize}
 
	We based $\mathscr{T}^\rho_*$ using the convention of Subsection~\ref{AlgPre}.

	First, the compatibility of the volume forms on the components of $\mathscr{T}^\rho_*$ implies $\mathrm{Tor}(\mathscr{T}^\rho_*) = 1$.
	
	Next, to compare the Reidemeister torsion of $\mathscr{C}^\rho_*$ with the one of $\mathscr{T}^\rho_*$ we use the commutative diagram 
$$ \xymatrix@-.9pc{
      & 0 \ar[d] & 0 \ar[d] & 0 \ar[d] & \\
\mathscr{S}_* = 0 \ar[r] & \su \ar[r] \ar[d] & \su \oplus \su \ar[r] \ar[d] & \su \ar[r] \ar[d] & 0 \\
\mathscr{C}^\rho_* = 0 \ar[r] & \mathrm{Der}_\rho (M_K) \ar[r] \ar[d] & \mathrm{Der}^{\mathcal{T}_1}_{\rho_1} (B_1) \oplus \mathrm{Der}^{\mathcal{T}_2}_{\rho_2} (B_2) \ar[r] \ar[d] & \mathrm{Der}^{\mathcal{S}}_{\rho_S} (S) \ar[r] \ar[d] & 0  \\
\mathscr{T}^\rho_* = 0 \ar[r] & \tangent{[\rho]}{\widehat{R}(M_K)} \ar[r] \ar[d] & \tangent{[\rho_1]}{\widehat{R}^{\mathcal{T}_1}(B_1)} \oplus \tangent{[\rho_2]}{\widehat{R}^{\mathcal{T}_2}(B_2)} \ar[r] \ar[d] & \tangent{[\rho_S]}{\widehat{R}^{\mathcal{S}}(S)} \ar[r] \ar[d] & 0 \\
& 0  & 0  & 0  &
}$$
Each row and each column of the previous diagram is an acyclic chain complex.

\begin{lemma}\label{L:TorC}
	The Reidemeister torsion of $\mathscr{C}^\rho_*$ satisfies $\mathrm{Tor}(\mathscr{C}^\rho_*) = -1$.
\end{lemma}

\begin{proof}

	The spaces $\mathscr{S}_0$ and $\mathscr{S}_2$ are endowed with the $3$-volume form $\ii \wedge \jj \wedge \kk$, and the space $\mathscr{S}_1 = \su \oplus \su$ is endowed with the product volume form. 
By construction the volume forms on the components of the short exact sequence $\xymatrix@1@-.6pc{0 \ar[r] & \mathscr{S}_i \ar[r] & \mathscr{C}^\rho_i \ar[r] & \mathscr{T}^\rho_i \ar[r] & 0}$ are compatible, for $i=0, 1, 2$ (see~Subsection~\ref{SS:FormeVolNat}). 
	
	The Multiplicativity Lemma gives $\mathrm{Tor}(\mathscr{C}^\rho_*) = (-1)^\alpha \mathrm{Tor}(\mathscr{S}_*) \cdot \mathrm{Tor}(\mathscr{T}^\rho_*)$,
where $\alpha = \alpha(\mathscr{S}_*, \mathscr{T}^\rho_*) = 1 \in \ZZ/2\ZZ$.
	 It is easy to see that $\mathrm{Tor}(\mathscr{S}_*) = 1$, which gives $\mathrm{Tor}(\mathscr{C}^\rho_*) = -1.$
\end{proof}

\subsection{Second Step}
	 In this second part we compare the sign-determined Reide\-meis\-ter torsion of the chain complex 
\[
\mathscr{C}^{\rho}_* = \xymatrix@1@-.6pc{0 \ar[r] & \mathrm{Der}_\rho (M_K) \ar[r] & \mathrm{Der}^{\mathcal{T}_1}_{\rho_1} (B_1) \oplus \mathrm{Der}^{\mathcal{T}_2}_{\rho_2} (B_2) \ar[r] & \mathrm{Der}^{\mathcal{S}}_{\rho_S} (S) \ar[r] & 0}
\]
with the one of 
\[
\mathscr{B}^{\rho}_* = \xymatrix@1@-.5pc{0 \ar[r] & {\Der{\rho}{(M_K)}}_0 \ar[r]  & {\mathrm{Der}_{\rho_1}^{\mathcal{T}_1}(B_1)}_0 \oplus {\mathrm{Der}_{\rho_2}^{\mathcal{T}_2}(B_2)}_0 \ar[r] & {\mathrm{Der}_{\rho_S}^{\mathcal{S}}(S)}_0 \ar[r] & 0.}
\]
The chain complex $\mathscr{B}^{\rho}_*$ is acyclic if and only if ${\mathrm{Der}_\rho(M_K)}_0 = \mathrm{Inn}_\rho(M_K)$, i.e. if and only if  $\rho$ is $m$-regular (see Remark~\ref{remarkreg}). If not, the homology of $\mathscr{B}^{\rho}_*$ reduces to $H_*(\mathscr{B}^{\rho}_*) = H_0(\mathscr{B}^{\rho}_*) \cong H^2_\rho(M_K) \cong \IR$. We thus split the discussion into two cases to make it clearer. 

\medskip

\paragraph*{(1) \emph{$m$-regular case}} 
	We suppose that $\rho$ is $m$-regular.

	To compare the Reidemeister torsion of the acyclic chain complex $\mathscr{C}^\rho_*$ to the one of the acyclic chain complex $\mathscr{B}^\rho_{*}$ we use the following commutative diagram 
\[
\xymatrix@-.7pc{ & 0 \ar[d] & 0 \ar[d] & 0 \ar[d] & \\
\mathscr{B}^{\rho}_* = 0 \ar[r] & {\Der{\rho}{(M_K)}}_0 \ar[r] \ar[d] & {\mathrm{Der}_{\rho_1}^{\mathcal{T}_1}(B_1)}_0 \oplus {\mathrm{Der}_{\rho_2}^{\mathcal{T}_2}(B_2)}_0 \ar[r] \ar[d]& {\mathrm{Der}_{\rho_S}^{\mathcal{S}}(S)}_0 \ar[r] \ar[d] & 0 \\
\mathscr{C}^\rho_* = 0 \ar[r] & \Der{\rho}{(M_K)} \ar[r] \ar[d]^-\beta  & \mathrm{Der}_{\rho_1}^{\mathcal{T}_1}(B_1) \oplus \mathrm{Der}_{\rho_2}^{\mathcal{T}_2}(B_2) \ar[r] \ar[d]^-{\pi_1 \oplus \pi_2}& \mathrm{Der}_{\rho_S}^{\mathcal{S}}(S) \ar[r] \ar[d]^-{\pi_S} & 0 \\
\mathscr{R}_* = 0 \ar[r] & \IR \ar[r]^-\alpha \ar[d] & \IR \oplus \IR \ar[r]^-\sigma \ar[d] & \IR \ar[r] \ar[d] & 0\\
& 0 & 0 & 0 &}
\]
Here $\pi_i$ is defined by equation~(\ref{Pii}), $i=1,2$, $\pi_S$ by equation~(\ref{PiS}),
\[
\beta : \mathrm{Der}_\rho(M_K) \to \IR,\; d \mapsto \langle d(m), P^\rho(m)\rangle,\]
\[
\alpha : \IR \to \IR \oplus \IR, \; x \mapsto (x, -x) \text{ and } \sigma : \IR \oplus \IR \to \IR, \; (x, y) \mapsto x+y.
\]
Each row and each column of the previous diagram is an acyclic chain complex. 

	The commutativity of the previous diagram comes from the fact that all elements of the systems of generators $\mathcal{T}_i$, $i=1, 2$, and $\mathcal{S}$ are conjugate to the meridian of $K$.

	Recall that 
\begin{itemize}
  \item $\mathscr{B}^\rho_2 = {\Der{\rho}{(M_K)}}_0 = \mathrm{Inn}_\rho(M_K) \cong \su$ is endowed with the $3$-volume form induced by the one on $\su$,
  \item $\mathscr{B}^\rho_1 = {\mathrm{Der}_{\rho_1}^{\mathcal{T}_1}(B_1)}_0 \oplus {\mathrm{Der}_{\rho_2}^{\mathcal{T}_2}(B_2)}_0$ is endowed with the $4n$-volume form induced by the one on $\mathscr{C}^\rho_1 = \mathrm{Der}_{\rho_1}^{\mathcal{T}_1}(B_1) \oplus \mathrm{Der}_{\rho_2}^{\mathcal{T}_2}(B_2)$,
  \item $\mathscr{B}^\rho_0 = {\mathrm{Der}_{\rho_S}^{\mathcal{S}}(S)}_0$ is endowed with the $(4n-3)$-volume form induced by the one on $\mathscr{C}^\rho_0 = \mathrm{Der}_{\rho_S}^{\mathcal{S}}(S)$ (see Subsection~\ref{Tspace}).
\end{itemize}

	As $\rho$ is supposed to be $m$-regular there exists a unique vector $h^{(1)}_\rho(m) \in  H^1_\rho(M_K)$ which satisfies $f_m^\rho(h^{(1)}_\rho(m)) = \langle h^{(1)}_\rho(m).m, P^\rho(m) \rangle = 1$, see Subsection~\ref{RegularRep}. Set 
\begin{equation}\label{tau}
\tau = \omega^K_{[\rho]}\left(\varphi^{-1}_{[\rho]}(h^{(1)}_\rho(m))\right)
\end{equation}
 which is expected to be the $(Ad \circ \rho)$-twisted Reidemeister torsion of $M_K$ with respect to the basis $\left\{ h^{(1)}_\rho(m), h^{(2)}_\rho\right\}$ of $H^*_\rho(M_K)$. With this notation we compute the Reidemeister torsion of $\mathscr{B}_*^{\rho}$ in this case. 
 
\begin{lemma}\label{L:TorD0}
	 We have $\mathrm{Tor}(\mathscr{B}^{\rho}_*) = \tau$. 
\end{lemma}
\begin{proof}
	The volume forms on the components of the short exact sequence $$\xymatrix@1@-.5pc{0 \ar[r] & \mathscr{B}^{\rho}_i \ar[r] & \mathscr{C}^\rho_i \ar[r] & \mathscr{R}_i \ar[r] & 0}$$ are compatible for $i =  0, 1$. We easily show that the volume form on $\mathscr{R}_2 = \IR$ which is compatible with $\ii \wedge \jj \wedge \kk$ on $\mathscr{B}^{\rho}_2 = {\mathrm{Der}_\rho(M_K)}_0 = \mathrm{Inn}_\rho(M_K) \cong \su$ and with $w_\rho$ on $\mathscr{C}^\rho_2 = \mathrm{Der}_\rho(M_K)$ is the usual $1$-volume form on $\IR$ {corrected} by the factor $\tau^{-1}$.

	The Multiplicativity Lemma applies and gives $$\mathrm{Tor}(\mathscr{C}^\rho_*) = (-1)^\alpha \mathrm{Tor}(\mathscr{B}^{\rho}_*) \cdot \mathrm{Tor}(\mathscr{R}_*),$$ where $\alpha = \alpha(\mathscr{B}^{\rho}_*, \mathscr{R}_*) = 1 \in \ZZ/2\ZZ$. It is easy to see that $\mathrm{Tor}(\mathscr{R}_*) = \tau^{-1}$. Thus   $\mathrm{Tor}(\mathscr{C}^\rho_*) = - \tau^{-1}\cdot \mathrm{Tor}(\mathscr{B}^{\rho}_*)$. We conclude referring ourselves to Lemma~\ref{L:TorC}.
\end{proof}

\medskip

\paragraph*{(2) \emph{Non $m$-regular case}} 
	Now we suppose that $\rho$ is regular but not $m$-regular.
	
	To compare the Reidemeister torsion of $\mathscr{C}^\rho_*$ with the sign-determined one of the non acyclic based and homology based chain complex $\mathscr{B}^{\rho}_*$ we use the commutative diagram
\[
\xymatrix@-.7pc{ & 0 \ar[d] & 0 \ar[d] & 0 \ar[d] & \\
\mathscr{B}^{\rho}_* = 0 \ar[r] & \Der{\rho}{(M_K)} \ar[r] \ar[d] & {\mathrm{Der}_{\rho_1}^{\mathcal{T}_1}(B_1)}_0 \oplus {\mathrm{Der}_{\rho_2}^{\mathcal{T}_2}(B_2)}_0 \ar[r] \ar[d]& {\mathrm{Der}_{\rho_S}^{\mathcal{S}}(S)}_0 \ar[r] \ar[d] & 0 \\
\mathscr{C}^\rho_* = 0 \ar[r] & \Der{\rho}{(M_K)} \ar[r] \ar[d]  & \mathrm{Der}_{\rho_1}^{\mathcal{T}_1}(B_1) \oplus \mathrm{Der}_{\rho_2}^{\mathcal{T}_2}(B_2) \ar[r] \ar[d]^-{\pi_1 \oplus \pi_2}& \mathrm{Der}_{\rho_S}^{\mathcal{S}}(S) \ar[r] \ar[d]^-{\pi_S} & 0 \\
\mathscr{R}'_* = 0 \ar[r] & 0 \ar[r] \ar[d] & \IR \oplus \IR \ar[r]^-\sigma \ar[d] & \IR \ar[r] \ar[d] & 0\\
& 0 & 0 & 0 &}
\]
Here $\sigma : \IR \oplus \IR \to \IR$ is defined by $\sigma(x, y) = x+y$ and we recall that $\mathscr{B}^{\rho}_2 = \Der{\rho}{(M_K)} = {\Der{\rho}{(M_K)}}_0$ in this case.

	Let $\mathcal{H}_*$ denote the long exact sequence in homology associated to the exact sequence of chain complexes $\xymatrix@1@-.5pc{0 \ar[r] & \mathscr{B}^{\rho}_* \ar[r] & \mathscr{C}^\rho_* \ar[r] & \mathscr{R}'_* \ar[r] & 0}$. Here $\mathcal{H}_*$ reduces to
\begin{equation}\label{EQH}
\mathcal{H}_* = \xymatrix@-.5pc{0 \ar[r] & H_1(\mathscr{R}'_*) \ar[r]^-\cong & H_0(\mathscr{B}^\rho_*) \ar[r] & 0.}
\end{equation}

	We have $H_1(\mathscr{R}'_*) = \ker \sigma =\{(x, -x)\;|\; x \in \IR\}$. Choose $(1, -1)$ as generator for $H_1(\mathscr{R}'_*) = \ker \sigma$. This choice will not affect our computation as we will see in Claim~\ref{Claim24}. With this choice we have (for the same reason as in the previous case) $\mathrm{Tor}(\mathscr{R}'_*) = 1$ (because $|\mathscr{R}'_*| = 0 \in \ZZ/2\ZZ$). 

\begin{lemma}\label{TorD0}
	If $H_0(\mathscr{B}^{\rho}_*)$ is based with any generator, then the sign-determined Reidemeister torsion of $\mathscr{B}^{\rho}_*$ satisfies: $\mathrm{Tor}(\mathscr{B}^{\rho}_*) = (\mathrm{tor}(\mathcal{H}_*))^{-1}$.
\end{lemma}
\begin{proof}
	As in the proof of Lemma~\ref{L:TorD0} we verify that the volume forms on the components of the exact sequence $$\xymatrix@1@-.5pc{0 \ar[r] & \mathscr{B}^{\rho}_i \ar[r] & \mathscr{C}^\rho_i \ar[r] & \mathscr{R}'_i \ar[r] & 0}$$ are compatible for $i=0, 1, 2$. 
	
	The Multiplicativity Lemma applies and gives
$$\mathrm{Tor}(\mathscr{C}^\rho_*) = (-1)^{\varepsilon + \alpha} \mathrm{Tor}(\mathscr{B}^{\rho}_*) \cdot \mathrm{Tor}(\mathscr{R}'_*) \cdot \mathrm{tor}(\mathcal{H}_*),$$ 
where $\varepsilon = \varepsilon(\mathscr{B}^{\rho}_*, \mathscr{C}^\rho_*, \mathscr{R}'_*) = 1 \in \ZZ/2\ZZ$ and $\alpha = \alpha(\mathscr{B}^{\rho}_*, \mathscr{R}'_*) = 0 \in \ZZ/2\ZZ$.
	
	Finally, Lemma~\ref{L:TorC} provides $\mathrm{Tor}(\mathscr{B}^{\rho}_*) = (\mathrm{tor}(\mathcal{H}_*))^{-1}.$
\end{proof}

\subsection{Third Step}
	This third and last step of the proof of Lemma~\ref{P:TorD} is the more difficult one. We compare the sign-determined Reidemeister torsion of 
\[
\mathscr{B}^{\rho}_* = \xymatrix@1@-.5pc{0 \ar[r] & {\Der{\rho}{(M_K)}}_0 \ar[r]  & {\mathrm{Der}_{\rho_1}^{\mathcal{T}_1}(B_1)}_0 \oplus {\mathrm{Der}_{\rho_2}^{\mathcal{T}_2}(B_2)}_0 \ar[r] & {\mathrm{Der}_{\rho_S}^{\mathcal{S}}(S)}_0 \ar[r] & 0.}
\]
with the one of the based and homology based chain complex (see Subsection~\ref{AuxDerCW})
$$\mathscr{D}^{\rho}_* = \xymatrix@1@-.5pc{0 \ar[r] & \Der{\rho}{(M_K)} \ar[r]^-{p^*}  & \Der{\rho_1}{(B_1)} \oplus \Der{\rho_2}{(B_2)} \ar[r]^-{\kappa^*}& \Der{\rho_S}{(S)} \ar[r] & 0}$$

	With the notation and results obtained in the last two steps we can achieve the computation of $\mathrm{Tor}(\mathscr{D}^\rho_*)$.
	
\begin{proof}[Proof of Lemma~\ref{P:TorD}]
	It requires several steps. 
	We must repeat the splitting used in the second step to use the computations we made for $\mathrm{Tor}(\mathscr{B}^{\rho}_*)$.

\medskip

\paragraph*{(1) \emph{$m$-regular case}}
	We suppose that $\rho$ is $m$-regular.
	
	Consider the linear maps
\begin{align*}
	\beta &: \mathrm{Der}_\rho(M_K) \to \IR, \; d \mapsto \langle d(m), P^\rho(m)\rangle,\\
	\beta_i &: \mathrm{Der}_{\rho_i}(B_i) \to \IR^n, \; d \mapsto \left({\langle d(t^{(i)}_1), P^{\rho_i}(t^{(i)}_1)\rangle, \ldots, \langle d(t^{(i)}_n), P^{\rho_i}(t^{(i)}_n)\rangle}\right) \; (i=1, 2),\\
	\beta_S &: \mathrm{Der}_{\rho_S}(S) \to \IR, \; d \mapsto \left({\langle d(s_1), P^{\rho_S}(s_1)\rangle, \ldots, \langle d(s_{2n}), P^{\rho_S}(s_{2n})\rangle}\right).
\end{align*}
We can easily see that these three maps are surjective and that
\[
{\mathrm{Der}_{\rho}(M_K)}_0 = \ker \beta, \; {\mathrm{Der}^{\mathcal{T}_i}_{\rho_i}(B_i)}_0 = \ker \beta_i, \; {\mathrm{Der}^{\mathcal{S}}_{\rho_S}(S)}_0 = \ker \beta_S.
\]
If $\psi : H^2_\rho(M_K) \to \IR$ denotes the isomorphism defined by $\psi(h^{(2)}_\rho) = 1$, then the following diagram is commutative:
$$ \xymatrix@-1.4pc{
      & 0 \ar[d] & 0 \ar[d] & 0 \ar[d] & 0 \ar[d] & \\
\overline{\mathscr{B}}^{\rho}_* = 0 \ar[r] & {\mathrm{Der}_\rho (M_K)}_0 \ar[r] \ar[d] & {\mathrm{Der}^{\mathcal{T}_1}_{\rho_1}(B_1)}_0 \oplus {\mathrm{Der}^{\mathcal{T}_2}_{\rho_2}(B_2)}_0 \ar[r] \ar[d] & {\mathrm{Der}^{\mathcal{S}}_{\rho_S}(S)}_0 \ar[r] \ar[d] & 0 \ar[r] \ar[d] & 0 \\
\overline{\mathscr{D}}^{\rho}_* = 0 \ar[r] & \Der{\rho}{(M_K)} \ar[r] \ar[d]^-{\beta} & \Der{\rho_1}{(B_1)} \oplus \Der{\rho_2}{(B_2)} \ar[r] \ar[d]^-{\beta_1 \oplus \beta_2} & \Der{\rho_S}{(S)} \ar[r] \ar[d]^-{\beta_S} & H^2_\rho(M_K) \ar[r] \ar[d]^-\psi & 0 \\
\overline{\mathscr{R}}_* = 0 \ar[r] & \IR \ar[r] \ar[d] & \IR^{n} \oplus \IR^n   \ar[r] \ar[d] &  \IR^{2n} \ar[r] \ar[d] & \IR \ar[r] \ar[d] & 0  \\
& 0  & 0  & 0  & 0 &
}$$

Each row of the previous diagram is an acyclic  chain complex. The volume forms induced on the components of $\overline{\mathscr{R}}_*$ which are compatible alongcolumns are the following ones: 
\begin{itemize}
  \item $\overline{\mathscr{R}}_3 = \IR$ is endowed with the usual $1$-volume form,
  \item $\overline{\mathscr{R}}_2 = \IR^{n} \oplus \IR^n$ is endowed with the usual product volume form on $\IR^n \oplus \IR^n$, 
  \item $\overline{\mathscr{R}}_1 = \IR^{2n}$ is endowed with the {opposite} of the usual product volume form,
  \item $\overline{\mathscr{R}}_0 = \IR$ is endowed with the usual $1$-volume form {corrected} by the factor $\tau^{-1}$ (see equation~(\ref{tau})). 
\end{itemize}

	First, observe that $\mathrm{Tor}(\mathscr{D}^\rho_*) = - {(\mathrm{Tor}(\overline{\mathscr{D}}^{\rho}_*))}^{-1}$ and $\mathrm{Tor}(\overline{\mathscr{B}}^{\rho}_*) = {(\mathrm{Tor}(\mathscr{B}^{\rho}_*))}^{-1}$. Moreover, Lemma~\ref{L:TorD0} implies $\mathrm{Tor}(\mathscr{B}^{\rho}_*) =  \tau$, thus $$\mathrm{Tor}(\overline{\mathscr{B}}^{\rho}_*) = \tau^{-1}.$$
	
	Next, the Multiplicativity Lemma gives
$$\mathrm{Tor}(\overline{\mathscr{D}}^{\rho}_*) =  - \tau^{-1} \cdot \mathrm{Tor}(\overline{\mathscr{R}}_*).$$
 It remains to compute $\mathrm{Tor}(\overline{\mathscr{R}}_*)$.
	For this purpose, let $\mathscr{M}_*$ be the chain complex equal to $\overline{\mathscr{R}}_*$ but differently based. We suppose that each component of $\mathscr{M}_*$  is endowed with the usual product volume form on $\IR^p$. If we apply the basis change formula~(\ref{EQ:changementdebase}) we obtain $$\mathrm{Tor}(\overline{\mathscr{R}}_*) = -\tau \cdot \mathrm{Tor}(\mathscr{M}_*).$$ The aim of the following claim is to compute the sign-determined Reidemeister torsion of $\mathscr{M}_*$.
\begin{claim}\label{Assertion}
 We have 
\begin{equation}\label{EQClaim}
\mathrm{Tor}(\mathscr{M}_*)= - \mathrm{sgn}(\mathrm{Tor}(\mathfrak{X}_*)).
\end{equation}
\end{claim}
Here we recall that $\mathfrak{X}_*$ is the chain complex~(\ref{NoeudReel}).
\begin{proof}[Proof of Claim~\ref{Assertion}]
	Several steps are required to accomplish it. 

\noindent (\textbf{a}) We begin by interpreting the acyclic  chain complex $\mathscr{M}_*$ as a Mayer-Vietoris sequence as follows. Consider a $3$-dimensional ball $D^3$ which contain the knot $K$. We write $\widehat{M}_K=M_K \cap D^3$, $\widehat{B}_1=B_1 \cap D^3$ and $\widehat{B}_2=B_2 \cap D^3$. 
	
	We can easily see that $H^1(\widehat{M}_K; \IR) \cong H^1(M_K; \IR) \cong \IR$, $H^1(\widehat{B}_i; \IR) \cong \IR^n$, $i=1, 2$. Using the long exact sequence in cohomology associated to the pair  $(\bord M_K, \widehat{M}_K)$, we further observe that $H^2(\widehat{M}_K; \IR) \cong H^2(\bord M_K; \IR) \cong \IR$.
	Thus the Mayer-Vietoris sequence with real coefficients corresponding to the splitting $\widehat{M}_K = \widehat{B}_1 \cup_{D^2\setminus K} \widehat{B}_2$  reduces to
\begin{equation}\label{MVchapeau}
\xymatrix@1@-.9pc{0 \ar[r] & H^1(\widehat{M}_K; \IR) \ar[r]^-{d^*} & H^1(\widehat{B}_1; \IR) \oplus H^1(\widehat{B}_2; \IR) \ar[r]^-{s^*} & H^1(D^2 \setminus K; \IR) \ar[r]^-\bord & H^2(\widehat{M}_K; \IR) \ar[r] & 0}.
\end{equation}

	To prove that the Mayer-Vietoris sequence~(\ref{MVchapeau}) is equal to $\mathscr{M}_*$ we introduce some more notation. The boundary of $B_i$, $i=1, 2$, is the union of $n$ disjoint cylinders denoted $T^{(i)}_j$. The boundary of $S$ is the disjoint union of $2n$ circles denoted $S^1_j$. The basis of $T^{(2)}_i$ are the circles $S^1_{2i-1}$ and $S^1_{2i}$ and the basis of $T^{(1)}_i$ are the circles $S^1_{j_{i, 1}}$ and $S^1_{j_{i, 2}}$, where  $j_{i,1}$ and $j_{i,2}$ are the two integers whose images by $\zeta$ are in the $i$th-handle (see Fig.~\ref{FigPlatGen} and Fig.~\ref{Fig:Tubes}).
	
	Let us furthermore consider the canonical inclusions: $i : \bord M_K \hookrightarrow M_K$, $j_i : S^1_i  \hookrightarrow S$, $k^{(1)}_i : T^{(1)}_i  \hookrightarrow \widehat{B}_1$, $k^{(2)}_i : T^{(2)}_i  \hookrightarrow \widehat{B}_2$, $l_i : S^1_i \hookrightarrow D^2 \setminus K$ and $\iota : \widehat{M}_K \hookrightarrow \bord M_K$. Write the Mayer-Vietoris sequences for cohomology associated to $M_K = B_1 \cup_S B_2$, to $\bord M_K = \bord B_1 \cup_{\bord S} \bord B_2$ and to $\widehat{M}_K = \widehat{B}_1 \cup_{D^2\setminus K} \widehat{B}_2$. Taking into account the naturality of the constructions and the fact that $H^1_{\rho_S}(S_i^1) = \IR P^{\rho_S}(s_i)$ and $H^1_\rho(\bord M_K) = \IR P^\rho(m)$, these sequences combine to yield the commutative diagram:
\[
\xymatrix@-.97pc{ H^1_{\rho_1}(B_1) \oplus H^1_{\rho_2}(B_2) \ar[r] \ar[d]_-{f_1 \oplus f_2} & H^1_{\rho_S}(S) \ar[r]^-\delta \ar[d]_-{\underset{i=1}{\overset{2n}{\oplus}} j^*_i}& H^2_{\rho}(M_K) \ar[d]^-{i^*} \ar[r] & 0\\
\displaystyle{\bigoplus_{i=1}^n H^1_{\rho_1}(T^{(1)}_i) \oplus \bigoplus_{i=1}^n H^1_{\rho_2}(T^{(2)}_i)} \ar[r] \ar[d]_-{r_1 \oplus r_2} & \displaystyle{\bigoplus_{i=1}^{2n} H^1_{\rho_S}(S^1_i)} \ar[r]^-{\Delta_1}\ar[d]_-{\underset{i=1}{\overset{2n}{\oplus}} P^{\rho_S}(s_i) \cup \, \cdot} & H^2_\rho(\bord M_K) \ar[d]^-{P^\rho(m) \cup \, \cdot} \ar[r] & 0\\
\displaystyle{\bigoplus_{i=1}^n H^1(T^{(1)}_i; \IR) \oplus \bigoplus_{i=1}^n H^1(T^{(2)}_i; \IR)} \ar[r]  &  \displaystyle{\bigoplus_{i=1}^{2n} H^1(S^1_i; \IR)} \ar[r]^-{\Delta_2} & H^2(\bord M_K; \IR) \ar[r] & 0\\
H^1(\widehat{B}_1) \oplus H^1(\widehat{B}_2) \ar[r] \ar[u]^-{\iota_1 \oplus \iota_2} &
H^1(D^2\setminus K; \IR) \ar[r]^-{\bord} \ar[u]^-{\underset{i=1}{\overset{2n}{\oplus}} l_i^*} & H^2(\widehat{M}_K; \IR) \ar[u]_-{\iota^*} \ar[r] & 0}
\]
In this diagram each map $\delta$, $\Delta_1$ and $\Delta_2$, and $\bord$ is the connecting operator of the Mayer-Vietoris sequence associated to $M_K = B_1 \cup_S B_2$, $\bord M_K = \bord B_1 \cup_{\bord S} B_2$ and $\widehat{M}_K = \widehat{B}_1 \cup_{D^2\setminus K} \widehat{B}_2$ respectively. We also use the following notation: 
\begin{align*}
f_1 = \bigoplus_{i=1}^n {(k^{(1)}_i)}^*,\quad & f_2 = \bigoplus_{i=1}^n {(k^{(2)}_i)}^*, \\ 
r_1 = \bigoplus_{i=1}^n P^{\rho_1}(t^{(1)}_i) \cup \cdot,\quad & r_2 = \bigoplus_{i=1}^n P^{\rho_2}(t^{(2)}_i) \cup \cdot, \\ 
\iota_1 = \bigoplus_{i=1}^n {(k^{(1)}_i)}^*,\quad & \iota_2 = \bigoplus_{i=1}^n {(k^{(2)}_i)}^*.
\end{align*}
and the following identifications $H^2(\bord M_K; \IR)  \cong \IR,$
$$\bigoplus_{i=1}^n H^1(T^{(1)}_i; \IR) \; \oplus \; \bigoplus_{i=1}^n H^1(T^{(2)}_i; \IR) \cong \IR^n \oplus \IR^n \text{ and } \bigoplus_{i=1}^{2n}H^1(S^1_i; \IR)  \cong \IR^{2n}.$$ 
All this suffices to establish the identification of $\mathscr{M}_*$ with the chain complex~(\ref{MVchapeau}).\\

\begin{figure}[!hbt]
\begin{center}
\begin{pspicture}(3, 5.5)
\psarc(.5, 5){.75}{0}{180}
\psarc(.5, 5){.25}{0}{180}
\psline(-.25,5)(-.25,4.75)\psline(.25,5)(.25,4.75)
\psline(.75,5)(.75,4.75)\psline(1.25,5)(1.25,4.75)
\pscurve(-.25,4.75)(-.2,4.7)(0,4.65)(.2,4.7)(.25,4.75)
\pscurve[linestyle=dotted](-.25,4.75)(-.2,4.8)(0,4.85)(.2,4.8)(.25,4.75)
\pscurve(.75,4.75)(.8,4.7)(1,4.65)(1.2,4.7)(1.25,4.75)
\pscurve[linestyle=dotted](.75,4.75)(.8,4.8)(1,4.85)(1.2,4.8)(1.25,4.75)
\psarc(3, 5){.75}{0}{180}
\psarc(3, 5){.25}{0}{180}
\psline(3.75,5)(3.75,4.75)\psline(2.25,5)(2.25,4.75)
\psline(3.25,5)(3.25,4.75)\psline(2.75,5)(2.75,4.75)
\pscurve(2.25,4.75)(2.3,4.7)(2.5,4.65)(2.7,4.7)(2.75,4.75)
\pscurve[linestyle=dotted](2.25,4.75)(2.3,4.8)(2.5,4.85)(2.7,4.8)(2.75,4.75)
\pscurve(3.25,4.75)(3.3,4.7)(3.5,4.65)(3.7,4.7)(3.75,4.75)
\pscurve[linestyle=dotted](3.25,4.75)(3.3,4.8)(3.5,4.85)(3.7,4.8)(3.75,4.75)
\pscurve(-.25,4.5)(-.2,4.45)(0,4.4)(.2,4.45)(.25,4.5)
\pscurve(-.25,4.5)(-.2,4.55)(0,4.6)(.2,4.55)(.25,4.5)
\pscurve(.75,4.5)(.8,4.45)(1,4.4)(1.2,4.45)(1.25,4.5)
\pscurve(.75,4.5)(.8,4.55)(1,4.6)(1.2,4.55)(1.25,4.5)
\pscurve(2.25,4.5)(2.3,4.45)(2.5,4.4)(2.7,4.45)(2.75,4.5)
\pscurve(2.25,4.5)(2.3,4.55)(2.5,4.6)(2.7,4.55)(2.75,4.5)
\pscurve(3.25,4.5)(3.3,4.45)(3.5,4.4)(3.7,4.45)(3.75,4.5)
\pscurve(3.25,4.5)(3.3,4.55)(3.5,4.6)(3.7,4.55)(3.75,4.5)
\pscurve(.75,4.5)(.8,4.25)(1.8,3.6)(2.25,3.25)(2,3)(1.8,2.9)(.5,2.1)(-.2,1.5)(-.25,1)
\pscurve(1.25,4.5)(1.3,4.3)(1.75,4)(2.25,3.7)(2.75,3.25)(2.25,2.8)(1,2)(.3,1.5)(.25,1)
\psline(-.25,1)(-.25,.5)\psline(.25,1)(.25,.5)
\psarc(.5, .5){.75}{180}{0}
\psarc(.5, .5){.25}{180}{0}
\psline[linewidth=.1, linecolor=white](1.25,2.55)(1.55,2.8)
\psline[linewidth=.1, linecolor=white](1.5,2.275)(1.9,2.55)
\psline(.75,.75)(.75,.5)\psline(1.25,.75)(1.25,.5)
\pscurve(.75,.75)(.8,1.25)(1.5,1.8)(1.6,2.2)(.75,3.25)(1.45,3.8)
\pscurve(1.25,.75)(1.3,1.25)(2,1.85)(2.1,2.2)(1.25,3.15)(1.7,3.65)
\pscurve(2,3.85)(2.7,4.25)(2.75,4.5)
\pscurve(1.75,4)(2.2,4.25)(2.25,4.5)
\pscurve(-.25,4.5)(0,3)(.6,2.15)\pscurve(.25,4.5)(.5,3)(.95,2.4)
\psline[linewidth=.1, linecolor=white](1.2,1.55)(1.55,1.85)
\psline[linewidth=.1, linecolor=white](1.4,1.4)(1.75,1.65)
\pscurve(.8,1.9)(1.25,1.5)(2.2,1)(2.25,.75)
\pscurve(1.15,2.075)(1.75,1.65)(2.7,1.1)(2.75,.75)
\psline(3.25,4.5)(3.25,.5)\psline(3.75,4.5)(3.75,.5)
\psline(2.75,.75)(2.75,.5)\psline(2.25,.75)(2.25,.5)
\psarc(3, .5){.75}{180}{0}
\psarc(3, .5){.25}{180}{0}
\uput{0}[0](-1,5.5){ { $T_1^{(2)}$}}
\uput{0}[0](3.75,5.5){ {$T_2^{(2)}$}}
\uput{0}[0](-1,0){ {$T_1^{(1)}$}}
\uput{0}[0](3.75,0){ {$T_2^{(1)}$}}
\uput{0}[0](-.75,4.65){ {\small $S^1_1$}}
\uput{0}[0](1.2,4.65){ {\small $S^1_2$}}
\uput{0}[0](1.7,4.65){ {\small $S^1_3$}}
\uput{0}[0](3.75,4.65){ {\small $S^1_4$}}
\uput{0}[0](4.5,5){ { $\bord B_2 = T_1^{(2)} \cup T_2^{(2)}$}}
\uput{0}[0](4.5,1){ { $\bord B_1 = T_1^{(1)} \cup T_2^{(1)}$}}
\uput{0}[0](-5,3){ { $S^1_{j_{1,1}} = S^1_2, \; S^1_{j_{1,2}} = S^1_3$}}
\uput{0}[0](-5,2.5){ { $S^1_{j_{2,1}} = S^1_1, \; S^1_{j_{2,2}} = S^1_4$}}
\psarc[linewidth=.005, linecolor=black, linestyle=dashed](.5, 5){.5}{0}{180}
\psline[linewidth=.005, linecolor=black, linestyle=dashed](0,5)(0,4.75)\psline[linewidth=.005, linecolor=black, linestyle=dashed](1,5)(1,4.75)
\pscurve[linewidth=.005, linecolor=black, linestyle=dashed](0,4.5)(.25,3)(.775,2.3)
\pscurve[linewidth=.005, linecolor=black, linestyle=dashed](1.05,1.95)(1.55,1.55)(2.45,1)(2.5,.75)
\psline[linewidth=.005, linecolor=black, linestyle=dashed](2.5,.75)(2.5,.5)
\psarc[linewidth=.005, linecolor=black, linestyle=dashed](3, .5){.5}{180}{0}
\psline[linewidth=.005, linecolor=black, linestyle=dashed](3.5,4.5)(3.5,.5)
\psarc[linewidth=.005, linecolor=black, linestyle=dashed](3, 5){.5}{0}{180}
\psline[linewidth=.005, linecolor=black, linestyle=dashed](3.5,5)(3.5,4.75)
\psline[linewidth=.005, linecolor=black, linestyle=dashed](2.5,5)(2.5,4.75)
\pscurve[linewidth=.005, linecolor=black, linestyle=dashed](1,4.5)(1.05,4.25)(1.95,3.7)(2.5,3.25)(2.25,3)(1.7,2.7)
\pscurve[linewidth=.005, linecolor=black, linestyle=dashed](1.4,2.5)(1,2.2)(.1,1.5)(0,1)
\psline[linewidth=.005, linecolor=black, linestyle=dashed](0,1)(0,.5)
\psarc[linewidth=.005, linecolor=black, linestyle=dashed](.5, .5){.5}{180}{0}
\pscurve[linewidth=.005, linecolor=black, linestyle=dashed](1,.5)(1,1)(1.1,1.2)(1.25,1.4)
\pscurve[linewidth=.005, linecolor=black, linestyle=dashed](1.6,1.75)(1.85,2.1)(1.55,2.5)(1.05,3.3)(1.55,3.7)
\pscurve[linewidth=.005, linecolor=black, linestyle=dashed](1.95,3.95)(2.45,4.3)(2.5,4.5)
\uput{0}[0](1.4,3.25){ { $K$}}
\end{pspicture}
\caption{Boundary of the exterior of the figure eight knot and the tubes $T^{(i)}_j$}
\label{Fig:Tubes}
\end{center}
\end{figure}
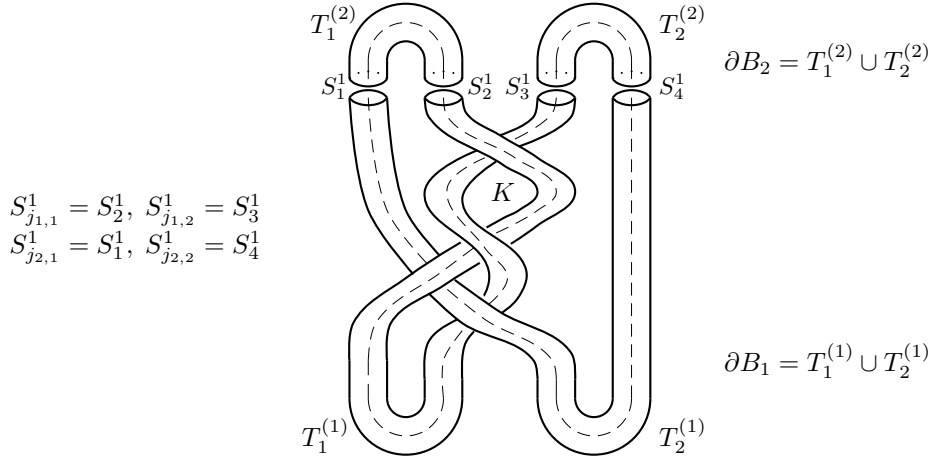

\noindent (\textbf{b}) Having the notation of part (\textbf{a}) in mind we are ready to compute $\mathrm{Tor}(\mathscr{M}_*)$. In this part we explicitly describe the homomorphisms $d^*$ and $s^*$ of the chain complex~(\ref{MVchapeau}).

Consider the following matrix corresponding to $d^*$
\[
P_{d^*} = \begin{pmatrix}
1 & 0 & \cdots & 0 & 0 & \cdots &0 \\ 
1 & 1 &  \cdots & 0 & 0&\cdots &0 \\ 
\vdots & \vdots & \ddots & \vdots & \vdots & & \vdots\\
1 & 0 & & 1 & 0& \cdots & 0\\
-1 & 0 & \cdots&0 & 1&  \cdots &0\\
\vdots & \vdots& & \vdots&\vdots &\ddots& \vdots\\
-1 & 0 & \cdots& 0&0 &\cdots &1
\end{pmatrix}.
\] 
Let $A$ denote the matrix of the homomorphism $s^*$; $A$ is a presentation matrix of $H^2(\widehat{M}_K; \mathbb{Z}) \cong \mathbb{Z}$. Let $P_{s^*}$ be the matrix deduced from 
$$A'=\begin{pmatrix} 
\; {A} & \left| \; \begin{matrix} 0 \\ \vdots \\ 1 \end{matrix} \right.
\end{pmatrix}$$
by deleting its first column.

As we easily see, we have 
\begin{equation}\label{torsionM}
\mathrm{Tor}(\mathscr{M}_*) = \det P_{s^*}\cdot  (\det P_{d^*})^{-1}.
\end{equation}

	We give an example to make the notation and the proof clearer.
\begin{example*}
	For the trefoil knot (see Fig.~\ref{FigPlatGen}) we have
\[
A=\begin{pmatrix}
1  & 0 & 1 & 0 \\
-1 & 0 & 0 & -1  \\
0 & -1 & -1 & 0 \\
0 & 1 & 0 & 1
\end{pmatrix} \quad \text{and} \quad 
P_{s^*}=\begin{pmatrix}
 0 & 1 & 0 & 0\\
 0 & 0 & -1 & 0 \\
 -1 & -1 & 0  & 0\\
 1 & 0 & 1 & 1
\end{pmatrix}.
\]

	For the figure eight knot (see Fig.~\ref{Fig:Tubes}) we have
\[
A=\begin{pmatrix}
1  & 0 & 0 & 1 \\
-1 & 0 & -1 & 0  \\
0 & 1 & 1 & 0 \\
0 & -1 & 0 & -1
\end{pmatrix} \quad \text{and} \quad 
P_{s^*}=\begin{pmatrix}
 0 & 0 & 1 & 0\\
 0 & -1 & 0 & 0 \\
 1 & 1 & 0  & 0\\
 -1 & 0 & -1 & 1
\end{pmatrix}.
\]
\end{example*}

	We can obtain the matrix $A$ in the following way. 
	
	The $i$th column of $A$, with $1 \leqslant i \leqslant n$, is the vector with all entries of which are zero except the one of index $2i-1$ and $2i$ which is equal to $(-1)^{\varepsilon_{2i-1}}$ and $(-1)^{\varepsilon_{2i}}$ respectively. Here $\varepsilon_k \in \{\pm 1\}$ depends on the orientation of $K$ and is defined by the rule $i_*(s^{(2)}_{2k-1}) = (t^{(2)}_k)^{\varepsilon_k}$ (see Fig.~\ref{FigPlatGen}).

 The $(n+i)$th column of $A$, with $1 \leqslant i \leqslant n$, is the vector with all entries of which are zero except  the ones of index $j_{i,1}$ and $j_{i,2}$ which are equal to $(-1)^{\varepsilon_{j_{i,1}}}$ and $(-1)^{\varepsilon_{j_{i,2}}}$ respectively. Here $j_{i,1}$ and $j_{i,2}$ are the two integers whose image by $\zeta$ are in the $i$th-handle (see Fig.~\ref{FigPlatGen}).

	We moreover observe that each row of $A$ contains exactly two coefficients non zero, equal to $\pm 1$; and that each column contains exactly two coefficients non zero, opposed and equal to $\pm 1$. Thus, the sum of all rows's entries is zero.  Seeing that $A$ is a presentation matrix of $\ZZ$, the determinant of the matrix $A_{n,n}$, deduced from $A$ by deleting the $n$th row and the $n$th column, is equal to $\pm 1$. This remark will be used in the sequel.

\medskip

\noindent (\textbf{c}) In this last step we prove that $\mathrm{Tor}(\mathscr{M}_*) = - (\mathrm{Tor}(\mathfrak{X}_*))^{-1}$.

The chain complex~(\ref{NoeudReel}) denoted $\mathfrak{X}_* $ reduces to
\[
\xymatrix@1{0 \ar[r] & \IR \ar[r]^-{\delta_2} & \IR^{2n} \ar[r]^-{\delta_1} & \IR^{2n-1} \ar[r] & 0,}
\]
where $\delta_2 = 0$. We have $|\mathfrak{X}_*| = 1$ and $H_*(\mathfrak{X}_*)$ is based with $\{\lbrack \! \lbrack pt \rbrack \! \rbrack, m^*\}$, see Subsection~\ref{CohomOr}.
To compute $\mathrm{Tor}(\mathfrak{X}_*)$ in terms of $\mathrm{Tor}(\mathscr{M}_*)$, we associate to $\mathfrak{X}_*$ the acyclic chain complex
\[
\mathfrak{X}'_* = \xymatrix@1{0 \ar[r] & H^1(M_K; \IR) \ar[r]^-{\sigma} & \IR^{2n} \ar[r]^-{\delta_1} & \IR^{2n-1}\oplus \IR \ar[r]^-{pr_{2n}} & \IR \ar[r] & 0.}
\]
Here, $H^1(M_K; \IR)$ is based with $m^*$, where $\sigma(m^*) = (1, \ldots, 1, -1, \ldots, -1) \in \IR^{2n}$, and the matrix of $\delta_1$ is obtained from $A$ by replacing the last row by $(0, \ldots, 0)$. 

	Equation~(\ref{EQ:Decalage}) implies $\mathrm{tor}(\mathfrak{X}'_*) = -(\mathrm{Tor}(\mathfrak{X}_*))^{-1} = - \mathrm{Tor}(\mathfrak{X}_*) \in \{\pm 1\}$, see Remark~\ref{remsign}. With obvious notation we have
\[
\mathrm{tor}(\mathfrak{X}'_*) = [b^0/c^0]^{-1} \cdot [b^1\widetilde{b}^0/c^1] \cdot [b^2 \widetilde{b}^1/c^2]^{-1}\cdot [\widetilde{b}^2/c^3],
\]
where $[b^0/c^0] = [\widetilde{b}^2/c^3] = 1$, $[b^1\widetilde{b}^0/c^1] = \det A_{n, n} = \det P_{s^*}$ and \[
[b^2 \widetilde{b}^1/c^2] = \begin{vmatrix}
1 & 1 & \ldots & 0  \\
\vdots &\vdots & \ddots & \vdots \\
-1 & 0 & \ldots & 1 \\
-1 &  0 & \ldots & 0
\end{vmatrix} = -(-1)^{2n+1} = 1. \]
It follows $\mathrm{tor}(\mathfrak{X}'_*) = \mathrm{Tor}(\mathscr{M}_*)$, which complete the proof of Claim~\ref{Assertion} keeping in mind Remark~\ref{remsign}.
\end{proof}

Claim~\ref{Assertion} provides $\mathrm{Tor}(\mathscr{D}^{\rho}_*) = - (\mathrm{Tor}(\mathscr{M}_*))^{-1} = \mathrm{sgn}(\mathrm{Tor}(\mathfrak{X}_*))$ which achieve the proof of Lemma~\ref{P:TorD} in the $m$-regular case.
\medskip

\paragraph*{(2) \emph{Non $m$-regular case}}
	Now we suppose that $\rho$ is regular but not $m$-regular.
	
	The comparison of the Reidemeister torsion of the based and homology based chain complex ${\mathscr{B}^{\rho}_*}$ with the one of the homology based chain complex $\mathscr{D}^\rho_*$ uses the commutative diagram: 
$$ \xymatrix@-.9pc{
      & 0 \ar[d] & 0 \ar[d] & 0 \ar[d] &  \\
{\mathscr{B}^{\rho}_*} = 0 \ar[r] & \mathrm{Der}_\rho (M_K) \ar[r] \ar[d] & {\mathrm{Der}^{\mathcal{T}_1}_{\rho_1}(B_1)}_0 \oplus {\mathrm{Der}^{\mathcal{T}_2}_{\rho_2}(B_2)}_0 \ar[r] \ar[d] & {\mathrm{Der}^{\mathcal{S}}_{\rho_S}(S)}_0 \ar[r] \ar[d] & 0 \\
\mathscr{D}^\rho_* = 0 \ar[r] & \Der{\rho}{(M_K)} \ar[r] \ar[d] & \Der{\rho_1}{(B_1)} \oplus \Der{\rho_2}{(B_2)} \ar[r] \ar[d]^-{\beta_1 \oplus \beta_2} & \Der{\rho_S}{(S)} \ar[r] \ar[d]^-{\beta_S} & 0 \\
\mathscr{R}''_* = 0 \ar[r] & 0 \ar[r] \ar[d] & \IR^{n} \oplus \IR^n   \ar[r]^-s \ar[d] &  \IR^{2n} \ar[r] \ar[d] & 0  \\
& 0  & 0  & 0  &
}$$
Let $\mathcal{H}'_*$ denote the long exact sequence in homology associated to the previous diagram; $\mathcal{H}'_*$ reduces to
\begin{equation}\label{H'}
\mathcal{H}'_* = \xymatrix@1@-.7pc{0 \ar[r] & H_1(\mathscr{R}''_*) \ar[r]^-{\cong} & H_0({\mathscr{B}^{\rho}_*}) \ar[r] & H_0(\mathscr{D}^\rho_*) \ar[r]^-\cong & H_0(\mathscr{R}''_*) \ar[r] & 0.}
\end{equation}
Remark that $H_2(\mathscr{R}''_*) = \ker s = \{(x, \ldots, x, -x, \ldots, -x) \;|\; x \in \IR\} \cong \IR$ and choose $(1, \ldots, 1, -1, \ldots, -1)$ as generator.

\begin{claim}\label{Claim24}
If $H_0(\mathscr{D}^\rho_*) = H^2_\rho(M_K)$ is based with the generator $h^{(2)}_\rho$ (see equation~(\ref{EQ:Defh2})) and if $\mathcal{H}_*$ denotes the chain complex~(\ref{EQH}), then $\mathrm{tor}(\mathcal{H}'_*) = \mathrm{tor}(\mathcal{H}_*)$.
\end{claim}
\begin{proof}
	We begin by the following observation
$$H_1(\mathscr{R}''_*) = \ker s = \Delta^{(n)} \oplus \Delta^{(n)} (\ker \sigma) = \Delta^{(n)} \oplus \Delta^{(n)} (H_1(\mathscr{R}'_*)),$$
where $\Delta^{(n)} : \IR \ni x \mapsto (x, \ldots, x) \in \IR^n$. Thus the isomorphism $\Delta^{(n)} \oplus \Delta^{(n)} : \IR \oplus \IR \to \IR^n \oplus \IR^n$ maps the generator $(1, -1)$ of $H_1(\mathscr{R}'_*)$ to the generator $(1, \ldots, 1, -1, \ldots, -1)$ of $H_1(\mathscr{R}''_*)$.
	
	Elsewhere, the Reidemeister torsion of the  chain complex $\mathcal{H}'_*$ is the inverse of the product of two changes of basis: the one corresponding to the isomorphism $H_0(\mathscr{D}^\rho_*) \cong H_0(\mathscr{R}''_*)$ and the other to the isomorphism $H_1(\mathscr{R}''_*) \cong H_0(\mathscr{B}^{\rho}_*)$.
	Keeping in mind the choices of bases, the determinant of the transition matrix of the change of basis corresponding to the first isomorphism $H_0(\mathscr{D}^\rho_*) \cong H_0(\mathscr{R}''_*)$ is $+1$. We remark that $$\mathbf{x} = (x, \ldots, x, -x, \ldots, -x) =  \Delta^{(n)} \oplus \Delta^{(n)}(x, -x) \in \ker s$$ lifts to an element in $\mathrm{Der}^{\mathcal{T}_1}_{\rho_1}(B_1) \oplus \mathrm{Der}^{\mathcal{T}_2}_{\rho_2}(B_2)$. Thus the determinant of the transition matrix of the change of basis corresponding to the second isomorphism $H_1(\mathscr{R}''_*) \cong H_0(\mathscr{B}^{\rho}_*)$ is equal to the determinant of the transition matrix of the change of basis corresponding to the isomorphism $H_1(\mathscr{R}'_*) \cong H_0(\mathscr{B}^{\rho}_*)$. This is sufficient to prove the required equality.  
\end{proof}

	We can now finish the computation of $\mathrm{Tor}(\mathscr{D}^\rho_*)$.
	We establish the compatibility along each column of the previous diagram using the same technique as in the case (1) where $\rho$ is supposed $m$-regular. The Multiplicativity Lemma gives 
\[
\mathrm{Tor}(\mathscr{D}^\rho_*) = (-1)^{\varepsilon+\alpha}\mathrm{Tor}({\mathscr{B}^{\rho}_*}) \cdot \mathrm{Tor}(\mathscr{R}''_*) \cdot \mathrm{tor}(\mathcal{H}'_*)
\]
with $\varepsilon = \varepsilon({\mathscr{B}^{\rho}_*}, \mathscr{D}^\rho_*, \mathscr{R}''_*) = 0 \in \ZZ/2\ZZ$ and $\alpha = \alpha(\mathscr{B}^{\rho}_*, \mathscr{R}''_*) = 0 \in \ZZ/2\ZZ$.

Claim~\ref{Claim24} and Lemma~\ref{TorD0} respectively furnish $\mathrm{tor}(\mathcal{H}'_*) = \mathrm{tor}(\mathcal{H}_*)$ and $\mathrm{Tor}(\mathscr{B}^{\rho}_*) = (\mathrm{tor}(\mathcal{H}_*))^{-1}$. Thus $\mathrm{Tor}(\mathscr{D}^\rho_*) = \mathrm{Tor}(\mathscr{R}''_*)$. It just remains to compute $\mathrm{Tor}(\mathscr{R}''_*)$. Using Claim~\ref{Assertion} it is easy to see that $\mathrm{Tor}(\mathscr{R}''_*) = \mathrm{sgn}(\mathrm{Tor}(\mathfrak{X}_*))$, which completes the proof of Lemma~\ref{P:TorD}.
\end{proof}

\section{Proof of Theorem~\ref{T:Volume=Torsion} : Part 2}
\label{S:Volume=Torsion2}

	The second part of the proof of Theorem~\ref{T:Volume=Torsion} consists in the comparison of the Reidemeister torsion of the chain complex
$$\mathscr{D}^\rho_* = \xymatrix@1@-.7pc{0 \ar[r] & \Der{\rho}{(M_K)} \ar[r]^-{p^*} & \Der{\rho_1}{(B_1)} \oplus \Der{\rho_2}{(B_2)} \ar[r]^-{\kappa^*} & \Der{\rho_S}{(S)} \ar[r] & 0}$$
to the one of the chain complex~(\ref{twistX}) of the knot exterior
 $$\mathscr{X}^\rho_* = \xymatrix@1@-.5pc{
0 \ar[r] & \su \ar[r]^-{d^{\rho}_2} & \su^{2n} \ar[r]^-{d^{\rho}_1} & \su^{2n-1} \ar[r] & 0}.$$
 It splits into two steps.

\subsection{First Step}
	We compare the sign-determined Reidemeister torsion of the non acyclic homology based chain complex $\mathscr{X}^\rho_*$ with the one of $\mathscr{D}^\rho_*$ using the diagram
$$ \xymatrix@-.9pc{
    & 0 \ar[d] & 0 \ar[d] & 0 \ar[d] & \\
\mathscr{X}^\rho_* = 0 \ar[r] & \su \ar[r]^-{d^\rho_2} \ar[d] & \su^{2n} \ar[r]^-{d^\rho_1} \ar[d]^-{\varphi_1 \oplus (-\varphi_2)} & \su^{2n-1} \ar[r] \ar[d]^-{\varphi_S} & 0 \\
\mathscr{D}^\rho_* = 0 \ar[r] & \Der{\rho}{(M_K)} \ar[r]^-{p^*} \ar[d] & \Der{\rho_1}{(B_1)} \oplus \Der{\rho_2}{(B_2)} \ar[r]^-{\kappa^*} \ar[d] & \Der{\rho_S}{(S)} \ar[r] \ar[d] & 0 \\
\mathscr{E}_* = 0 \ar[r] & H^1_{\rho}(M_K) \ar[r] \ar[d] & 0 \ar[r] \ar[d] &  0  \ar[r] \ar[d] &  0  \\
& 0  & 0  & 0  &}$$
Each component of the rows is endowed with a reference volume form. These reference volume forms are compatible along each column (see Subsection~\ref{SS:Der}).
	With the notation above, we prove
\begin{lemma}\label{LemmaPart2}
	For all non zero vector $v$ in $H^1_\rho(M_K)$, we have 
\[
	\tau^K_\rho\left( \varphi^{-1}_{[\rho]}(v)\right) = -(\mathrm{Tor}(\mathscr{E}_*))^{-1}.
\]
\end{lemma}

	Before proving Lemma~\ref{LemmaPart2} we establish 
\begin{claim}\label{C:com}
	The previous diagram is commutative.
\end{claim}
\begin{proof}
	To establish the commutativity of the diagram we use an explicit description of the boundary operators of the  chain complex $\mathscr{X}^\rho_*$. Using the Wirtinger presentation~(\ref{EQ:Presentation}) of $G_K$ we recall that
$$d^\rho_2(x) = \left( (1-t^{(1)}_1) \circ x, \ldots, (1-t^{(1)}_n) \circ x , (1-t^{(2)}_1) \circ x, \ldots, (1-t^{(2)}_n) \circ x \right)$$
for all $x \in \su$, and
\[
d^\rho_1((x_j)_{1 \leqslant j\leqslant 2n}) = \left( \sum_{j=1}^n \fox{r_i}{t^{(1)}_j} \circ x_{2j-1} + \sum_{j=1}^n \fox{r_i}{t^{(2)}_j} \circ x_{2j}\right)_{1\leqslant i \leqslant 2n-1}
\]
for all $(x_j)_{1 \leqslant j\leqslant 2n} \in \su^{2n}$, where $r_i = \kappa_1(s_i) \kappa_2(s_i^{-1})$. 	
	Essentially the commutativity comes from the two following equalities
\[
\fox{r_i}{t^{(1)}_k} = \fox{\kappa_1(s_i)}{t^{(1)}_k} \; \text{and}\; \fox{r_i}{t^{(2)}_k} = -\fox{\kappa_2(s_i)}{t^{(2)}_k},
\]
\end{proof}

	Next, the twisted homology groups of $\mathscr{X}^\rho_*$, $\mathscr{D}^\rho_*$ and $\mathscr{E}_*$ respectively satisfies
	\[
H_0(\mathscr{X}^\rho_*) \cong H^2_\rho(M_K), \; H_1(\mathscr{X}^\rho_*) \cong H^1_\rho(M_K) \; \text{and} \; H_i(\mathscr{X}^\rho_*) = 0 \; \text{if}\; i\geqslant 2 \text{ (see \S~\ref{TwistedC})}.
\]
\[
H_0(\mathscr{D}^\rho_*) \cong H^2_\rho(M_K) \; \text{and} \; H_i(\mathscr{D}^\rho_*) = 0 \; \text{if}\; i\geqslant 1 \text{ (Fact~\ref{HofD})}.
\]
\[
H_2(\mathscr{E}_*) \cong H^1_\rho(M_K) \; \text{and} \; H_i(\mathscr{E}_*) = 0 \; \text{if}\; i \ne 2.
\]

	Now we turn to the proof of Lemma~\ref{LemmaPart2}.

\begin{proof}[Proof of Lemma~\ref{LemmaPart2}]
	Let $\mathcal{H}''_*$ denote the long exact sequence in homology as\-so\-cia\-ted to $\xymatrix@1@-.7pc{
0 \ar[r] & \mathscr{X}^\rho_* \ar[r] & \mathscr{D}^\rho_* \ar[r] & \mathscr{E}_* \ar[r] & 0.}$ The Multiplicativity Lemma gives
$$\mathrm{Tor}(\mathscr{D}^\rho_*) = (-1)^{\varepsilon + \alpha}\mathrm{Tor}(\mathscr{X}^\rho_*) \cdot \mathrm{Tor}(\mathscr{E}_*) \cdot \mathrm{Tor}(\mathcal{H}''_*),$$
where $\varepsilon = \varepsilon(\mathscr{X}^\rho_*, \mathscr{D}^\rho_*, \mathscr{E}_*) = 0 \in \ZZ/2\ZZ$, $\alpha = \alpha(\mathscr{X}^\rho_*, \mathscr{E}_*) = 1  \in \ZZ/2\ZZ$.

The sequence $\mathcal{H}''_*$ reduces to the two isomorphisms: $\xymatrix@1@-.6pc{H_2(\mathscr{E}_*) \ar[r]^-{\cong} & H_1(\mathscr{X}^\rho_*)}$ and $\xymatrix@1@-.6pc{H_0(\mathscr{X}^\rho_*) \ar[r]^-{\cong} & H_0(\mathscr{D}^\rho_*)}$. Thus $\mathrm{Tor}(\mathcal{H}''_*)$ is the inverse of the product of two changes of basis, the one corresponding to $H^2_\rho(M_K)$ and the other to $H^1_\rho(M_K)$. Because of the choice of the bases of the twisted cohomology groups in dimension $1$ and $2$, we have $\mathrm{tor}(\mathcal{H}''_*) = 1$. So $\mathrm{Tor}(\mathscr{D}^\rho_*) = - \mathrm{Tor}(\mathscr{X}^\rho_*) \cdot \mathrm{Tor}(\mathscr{E}_*)$. Finally, using Lemma~\ref{P:TorD}, we conclude that
\begin{equation}\label{EQ3}
\mathrm{sgn}(\mathrm{Tor}(\mathfrak{X}_*)) \cdot \mathrm{Tor}(\mathscr{X}^\rho_*) = - {(\mathrm{Tor}(\mathscr{E}_*))}^{-1}.
\end{equation}
\end{proof}

\subsection{Second Step}
	It remains to bring together all the computations of Rei\-de\-meis\-ter torsions we have done before. 

\begin{proof}[End of the proof of Theorem~\ref{T:Volume=Torsion}]	
Fix a non zero vector 
$v \in H^1_\rho(M_K)$; we will prove
$$\vol{K}_{[\rho]}\left(\varphi^{-1}_{[\rho]}(v)\right) = \tors{K}_{\rho}\left(\varphi^{-1}_{[\rho]}(v)\right).$$
	
	To this end, we must compute the sign-determined Reidemeister torsion of the based and homology based chain complex $\mathscr{E}_*$. Let $z$ denote the generator of $H^1_\rho(M_K)$ such that $\vol{K}_{[\rho]}(\varphi^{-1}_{[\rho]}(z))=1$. The basis of $H^1_\rho(M_K)$ formed by the single vector $v$ and the basis of $H^2_\rho(M_K)$ formed by the single vector $h^{(2)}_\rho$ are the bases of the twisted cohomology groups of $M_K$ in which we compute the sign-determined Reidemeister torsion $\mathrm{Tor}(\mathscr{X}^\rho_*)$.	We observe that $\mathscr{E}_2 = H^1_\rho(M_K)$ is based with $\{z\}$ and $H_2(\mathscr{E}_*) \cong H^1_\rho(M_K)$ is based with $\{v\}$. With our conventions, we thus have $\mathrm{Tor}(\mathscr{E}_*) = - [z / v]$ (because $|\mathscr{E}_*| = 1$). 
	
	Lemma~\ref{LemmaPart2} provides
\begin{equation}\label{torsionE}
\tors{K}_{\rho}\left(\varphi^{-1}_{[\rho]}(v)\right) = - {(\mathrm{Tor}(\mathscr{E}_*))}^{-1} = [v / z].
\end{equation}
So, by using formula~(\ref{torsionE}) and the definition of $z$ we get
 $$\vol{K}_{[\rho]}\left(\varphi^{-1}_{[\rho]}(v)\right)= [v / z] =  \tors{K}_{\rho}\left(\varphi^{-1}_{[\rho]}(v)\right),$$
which completes the proof of Theorem~\ref{T:Volume=Torsion}.
\end{proof}
\subsection*{Acknowledgements} 
The author wishes to express his gratitude to Mi\-chael Heu\-se\-ner, Joan Porti, Daniel Lines, Vladimir Turaev, Michel Boileau, Jean-Yves Le\-Dimet, Louis Funar, Rinat Kashaev and Pierre de la Harpe for hepful discussions related to this paper.

The author also would like to thanks the referee for his careful reading of the paper and for his remarks which contribute to make it clearer.%

%
%
%
\bibliographystyle{amsalpha}
\bibliography{ReferencesArticle}
\end{document}